\documentclass[12pt]{iopart}

\usepackage{iopams}

\eqnobysec

\let \al=\alpha
\let \be=\beta

\let \var=\varphi
\let \vare=\varepsilon

\let \th=\theta
\let \la=\lambda
\let \de=\delta

\let \p=\partial
\let \q=\quad

\let \med=\medskip
\let \smal=\smallskip

\def\R{{\rm I\kern
-1.6pt{\rm R}}}
\def\C{{\rm |\kern
-4.6pt{\rm C}}}
\def\N{{\rm I\kern
-4.0pt{\rm N}}}

\def\q{\quad}

\def\ter{\hfill \vrule width 5 pt height 7 pt depth - 2 pt\smallskip}

\newtheorem{thm}{Theorem}[section]
\newtheorem{cor}[thm]{Corollary}
\newtheorem{lem}[thm]{Lemma}
\newtheorem{rem}[thm]{Remark}
\begin{document}

\title[Positive travelling fronts for delayed systems]{Positive travelling fronts for reaction-diffusion systems with distributed delay}

\author{Teresa Faria$^1$ and Sergei Trofimchuk$^2$}

\address{$^1$ Departamento de Matem\'atica  and CMAF, Faculdade de
Ci\^encias, Universidade de Lisboa, Campo Grande, 1749-016 Lisboa,
Portugal}
\address{$^2$ Instituto de Matem\'atica y Fisica,
Universidad de Talca, Casilla 747, Talca, Chile}
\eads{\mailto{tfaria@ptmat.fc.ul.pt},
\mailto{trofimch@inst-mat.utalca.cl}}
\begin{abstract}We give sufficient conditions for the existence of
positive travelling wave solutions  for  multi-dimensional
autonomous reaction-diffusion systems with distributed delay. To
prove the existence of travelling waves, we give an abstract
formulation   of the equation for the wave profiles in some suitable
Banach spaces, and  apply  known results about the index of some
associated Fredholm operators. After a Liapunov-Schmidt reduction,
these waves are obtained via the Banach contraction principle, as
perturbations of a positive heteroclinic solution  for the
associated system  without diffusion, whose existence is proven
under some requirements. By a careful analysis of the exponential
decay of the  travelling wave profiles at $-\infty$,  their
positiveness  is deduced. The existence  of positive travelling
waves is important in terms of applications to biological models.
Our method applies to systems of delayed reaction-diffusion
equations  whose nonlinearities are not required to satisfy a
quasi-monotonicity  condition. Applications are given, and include
the delayed Fisher-KPP equation.

\end{abstract}

\maketitle

\section{Introduction}

For the last decades, there has been an increasing number of studies
in travelling wave fronts for delayed diffusion equations, and
several methods
 to prove their existence have been developed.

In this paper, we are concerned with the existence and positiveness
of travelling waves  connecting two equilibria, for a class of
$N$-dimensional systems of reaction-diffusion equations with
distributed delay in the reaction terms, of the form
\begin{equation}\label{e1.1}
{{\p u}\over {\p t}}(t,x)=\Delta u(t,x)+f(u_t(\cdot ,x)),\q t\in
\mathbb{R},\  x\in\mathbb{R}^p.
\end{equation}
Here,  $f:{\cal C}:=C([-\tau ,0];\mathbb{R}^N)\to \mathbb{R}^N$ is
continuous, ${\cal C}$ is equipped with the norm
$\|\var\|_\infty=\sup_{\th\in [-\tau,0]} |\var(\th)|$, for some
fixed norm $|\cdot|$ in $\mathbb{R}^N$, and $\tau >0$.  As usual,
$u_t(\cdot ,x)$ denotes the restriction of a solution $u(t,x)$ to
the time interval $[t-\tau,t]$, i.e., $u_t(\th,x)=u(t+\th,x)$ for
$-\tau\le \th \le 0, x\in \mathbb{R}^p$. For simplicity, we consider
all the diffusion coefficients equal to 1 in (\ref{e1.1}), but all
our results apply to the more general case of  the diffusion term
given by $D\Delta u(t,x)$, where $D=diag \, (d_1,\dots,d_N)$ with
$d_i>0$.

We are mostly interested in situations where (1.1) represents a
population dynamics model, or another biological model. Typically,
we want to obtain conditions for the existence of a travelling front
connecting  two steady-states, zero and a positive equilibrium
$K\in\mathbb{R}^N$.  Due to the biological interpretation of the
model, only non-negative solutions are meaningful, therefore we look
for {\it positive} travelling wave solutions, connecting 0 to $K$ as
$t$ goes from $-\infty$ to $\infty$.

With the  method presented here, such positive travelling  waves are
obtained  for large wave speeds, as perturbations of a positive
heteroclinic solution for the corresponding functional differential
equation (FDE)  without diffusion,
\begin{equation}\label{e1.2}
u'(t)=f(u_t),\q t\in \mathbb{R}
\end{equation}
(where $u_t\in {\cal C}$ denotes the function $u_t(\th)=u(t+\th)$
for $\th \in [-\tau,0]$), whose existence we shall prove under some
requirements on $f$. This idea is not original, and has been
exploited in the literature (see e.g. \cite{FHW, FT}). When compared
with \cite{FHW}, in the present paper the major novelty is that we
give conditions for the travelling waves to be {\it positive}. We
also note that \cite{FHW} considers delayed reaction-diffusion
equations with a global space interaction, a situation not
considered here, for the sake of simplicity. Our results can however
be extended easily, to take into account non-local effects.  On the
other hand, \cite{FT}  deals with {\it scalar} reaction-diffusion
equations with one single {\it discrete} delay of the form
\begin{equation}\label{e1.3}
{{\p u}\over {\p t}}(t,x)=d{{\p^2 u}\over {\p x^2}}(t,x)+f(u(t,x),
u(t-\tau,x)),
\end{equation}
where $f(u(t,x), u(t-\tau,x)) =- u(t,x)+g(u(t-\tau,x))$ and
$g:[0,\infty)\to [0,\infty)$ is $C^2$-smooth, $g(0)=0, g(K)=K$ for
some $K>0$, and $g'(0)>1$.  Assuming that the two equilibria $0$ and
$K$ are  hyperbolic, under   some further assumptions  the existence
of positive and in general non-monotone travelling wavefronts
connecting 0 to $K$ was established in \cite{FT}.

Recently, several  techniques have been developed to prove the
existence of  travelling wave fronts for delayed diffusion
equations. They are often based on the application of a fixed point
theorem in an adequate Banach space, which requires a {\it
quasi-monotonicity} condition, either for the original equation
(\ref{e1.3})  \cite{M,WZ}  or, more recently,  for some auxiliary
equations \cite{MA}. These methods are usually combined with a
monotonic iteration scheme, associated with the construction of a
pair of upper and lower solutions. See \cite{M,MA,WZ} and references
therein. We emphasize that our method applies to systems
(\ref{e1.1}) with non-monotone nonlinearities, in the sense that we
do not impose on $f$ any type of quasi-monotonicity condition, as
defined in \cite{HS,WZ}.

Before introducing our hypotheses, we set some standard notation.
For $d=(d_1,\dots ,d_N)\in\mathbb{R}^N$, we say that $d>0$
(respectively $d\ge 0$) if $d_i>0$ (respectively $d_i\ge 0$) for
$i=1,\dots,N$. In ${\cal C}$, we consider the  partial order
$\phi\ge \psi$ if and only if $\phi(\th)- \psi(\th)\ge 0$ for
$\th\in [-\tau,0]$; in a similar way, $\phi >\psi$ if  $\phi(\th)-
\psi(\th)> 0$ for $\th\in [-\tau,0]$. As usual, ${\cal C}_+$ denotes
the positive cone $C([-\tau,0];[0,\infty)^N)$.

\med For $f$, the following hypotheses will be considered:
\begin{description}
\item [{(H1)}] $f(0)=f(K)=0$, where $K$ is some positive vector;
\item [{(H2)}] (i) $f$ takes bounded sets of ${\cal C}$ into bounded sets of $\mathbb{R}^N$ and is $C^2$-smooth; furthermore,
(ii) for all $M>0$ there is $\be >0$ such that $f_i(\var)+\be
\var_i(0)\ge 0,$ $ i=1,\dots,N$, for all $\var\in {\cal C}$ with
$0\le \var\le M$;
\item [{(H3)}] for Eq. (\ref{e1.2}),  the equilibrium $u=K$
is locally asymptotically stable and
 globally attractive in the set of solutions of (1.2)  with initial conditions $\var\in {\cal C}_+,\var(0)> 0$;
\item  [{(H4)}] for Eq. (\ref{e1.2}), its linearized equation about
the equilibrium 0 has a real characteristic root $\la_0>0$, which is
simple and dominant (i.e., $\Re\, z<\la_0$ for all other
characteristic roots $z$);  moreover,
  there is  a characteristic   eigenvector ${\bf v}>0$ associated  with $\la_0$.
\end{description}
We summarize the main results in this paper as follows. In Section
2, we assume (H1)-(H4) and establish the existence of a positive
heteroclinic solution $u^*(t)$ to (\ref{e1.2}), with
$u^*(-\infty)=0, u^*(\infty)=K$, and asymptotic behaviour
$O(e^{\la_0 t})$ at $-\infty$. In Section 3, for large wave speeds
we prove the existence of travelling wave solutions for
(\ref{e1.1}), connecting 0 to $K$. The  profiles of these waves are
obtained as perturbations of $u^*(t)$ via  a contraction principle
argument. For this, we generalize the procedure in \cite{FHW}, and
give an abstract formulation of the wave profiles as solutions of an
operational equation, acting in suitable Banach spaces, which
incorporate a desirable exponential decay  $O(e^{\mu t})$  at
$-\infty$, $0\le \mu<\la_0$.  Some nice results of Hale and  Lin
\cite{HL} on the index of some associated Fredholm operators are
used, and a Liapunov-Schmidt reduction effected, to set up the right
framework for the application of a contraction principle. As
mentioned above, an existence result of travelling waves connecting
two hyperbolic equilibria was already obtained in \cite{FHW}, for a
class of reaction-diffusion equations with global response, but for
such waves no exponential decay at $-\infty$  was derived in
\cite{FHW}, nor their positiveness. By a careful analysis of the
behaviour of the wave profiles at $-\infty$, in Section 4 we prove
that there are {\it positive} travelling waves if the wave speed is
large enough, and explicitly give their asymptotic decay at
$-\infty$. Section 5 is dedicated to applications, which include the
Fisher-KPP equation with delay and a 2-dimensional chemostat model.
An important theorem on the asymptotic behaviour of solutions of
perturbed linear autonomous ordinary FDEs is given in the Appendix.
This result generalizes to the case of FDEs with distributed  delay
a  result by Mallet-Paret \cite{MP}, for FDEs with discrete
time-delays (or time-shifts), and is often used in Sections 2 and 4.

\section{Positive heteroclinic solution for Eq. (\ref{e1.2})}

In this section, we prove the existence  of a positive solution
$u^*(t)$ of the ordinary FDE (\ref{e1.2}) connecting the equilibrium
0 to the positive equilibrium $K$. We recall that a function $u(s)$
defined on a set $S$ and with values in $\mathbb{R}^N$ is said to be
positive if all its components $u_1(s),\dots ,u_N(s)$ are positive
functions on $S$.
\begin{thm}\label{T2.1} Assume (H1)-(H4). Then:\vskip 0cm
i) There exists a heteroclinic
solution $u^*(t),t\in \mathbb{R},$ of Eq. (\ref{e1.2}), with
$u^*(-\infty )=0, u^*(\infty)=K$;\vskip 0cm (ii) $u^*(t)$ is
positive, $t\in\R$;\vskip 0cm (iii) $u^*(t)=ce^{\la_0t} {\bf
v}+O(e^{(2\la_0-\vare)t})$ at $-\infty$, for some $c>0$ and each
fixed $\vare>0$. \end{thm}

\noindent {\it Proof}. (i) Consider the linearization of (1.2) about
$0$,
\begin{equation}\label{e2.1}
u'(t)=Lu_t,\q {\rm where}\q  L=Df(0),
\end{equation}
and its characteristic equation
\begin{equation}\label{e2.2}
\det \Delta_0 (\la )=0,\q {\rm where}\q \Delta_0 (\la)=L(e^{\la
\cdot}I)-\la I.
\end{equation}
Recall that $\la $ is a solution of (\ref{e2.2}) if and only if
$\la\in \sigma (A)$, where $A$ is the infinitesimal generator
associated with the semi-flow of (\ref{e2.1}).

Let $\la_0>0$ be the leading (simple) eigenvalue of (\ref{e2.1})
given in (H4), and ${\bf v}\in\R^N, {\bf v}>0$, such that
$\Delta_0(\la_0){\bf v}=0$. Choose $\gamma>0$ with $\gamma
<\la_0<2\gamma$ and such that the strip $\gamma \le Re\, \la <\la_0$
does not contain any root of (\ref{e2.2}). Define
$\chi_0(\th)=e^{\la_0\th}{\bf v},\th \in [-\tau,0]$, and decompose
the phase ${\cal C}$ as ${\cal C}=P\oplus Q$, where $P=<\chi_0>$ and
$Q$ is the complementary space given by the formal adjoint theory of
Hale \cite{HVL}. Then there are neighbourhoods $N_0,N_1$ of 0 in
$P,Q$, respectively, and a $C^1$ map $w:N_0\to N_1$ with
$w(0)=0,Dw(0)=0$ such that the local $\gamma$-unstable manifold of 0
for Eq. (\ref{e1.2}) is given by
$$W(0)=\{ \phi +w(\phi): \phi\in N_0\}.$$
Note that  $\var \in W(0)$ if and only if  there is a full
trajectory $u_t=u_t(\var)\, ( t\in\mathbb{R})$ of (\ref{e1.2}) with
$u_0=\var$, $u_t\in N_0+N_1$ for $t\le 0$ and $u(t)e^{-\gamma t}\to
0$ as $t\to -\infty$. See Krisztin et al. \cite{KWW}, Hale and Lunel
\cite[Sec. 10.1-10.2]{HVL}, and Diekmann et al. \cite[Sec. 8.4]{DW}.

We now argue as in \cite{OW}. Write $w(t)=(w_1(t),\dots,w_N(t)),
{\bf v}=({\bf v}_1,\dots ,{\bf v}_N)$. Since $Dw(0)=0$, then
$\lim_{\|\phi\|\to 0, \phi\in N_0}{{\|w(\phi)\|}\over {\|\phi\|}}=0$
and  we have $\lim_{c\to 0}{{\|w(c \chi_0)\|}\over {|c|}}=0$. Thus,
there is $c_0>0$ such that $|w_i(c\chi_0)|_\infty\le ce^{-\la_0
\tau}{\bf v}_i/2$ for $c\in (0,c_0], i=1,\dots,N$, which implies
that for $0<c\le c_0$ we have
\begin{equation}\label{e2.3}
\min_{-\tau\le \th \le 0} \Big (ce^{\la_0 \th}{\bf
v}_i+w_i(c\chi_0)(\th)\Big )\ge ce^{-\la_0\tau}{\bf v}_i/2>0,\q
i=1,\dots,N,
\end{equation}
and therefore  $c\chi_0+w(c\chi_0)\in W(0)\cap {\cal C}_+$ for all
$c\in (0,c_0]$. Fix e.g. $c=c_0$, denote
$\phi=c_0\chi_0+w(c_0\chi_0)$ and  consider the full trajectory
$u_t^*=u_t(\phi),\ t\in \mathbb{R}$. We have  $u_t^*\in W(0)$ for
$t\le 0$, hence $u_t^*$ has the form
$u_t^*=c(t)\chi_0+w(c(t)\chi_0)$. Since the map $t\mapsto u_t^*$ and
the canonical projection of ${\cal C}$ on $P$ are continuous, $c(t)$
is continuous as well, with $c(t)\to 0$ as $t\to -\infty$. This
implies that there is $T<0$ such that $c(t)\le c_0$ for $t< T$. On
the other hand, if $c(t_0)=0$ for some $t_0<T$, then $u_{t_0}^*=0$,
which is not possible. From (\ref{e2.3}) it follows that $u^*(t)>0$
for $t<T$. Now, from (H3)  we have $u^*(t)\to K$ as $t\to\infty$.
This means that $u^*(t)$ is a heteroclinic solution of (\ref{e1.2})
connecting the two equilibria $0,K$, with $u^*(t)$ positive on some
interval $(-\infty,T)$.

(ii) Choose $M>0$ such that $u_i^*(t)\le M, \ t\in \mathbb{R},\
i=1,\dots,N$. For the sake of contradiction, suppose there is $t\ge
T$ and $i\in\{ 1,\dots,N\}:=I $ with $u_i^*(t)\le 0$. Define
$t^*=\min \{ t\ge T: u_j^*(t)= 0$ for some $j\in I\}$ and take $i\in
I$ such that $u_i^*(t^*)=0$. For $M$ as above, let $\be$ be as in
(H2), i.e., $f_j(\var)+\be \var_j(0)\ge 0,$ for  $j\in I$ and  $0\le
\var\le M$. Writing $u_i^*(t)$ in integral form,
\begin{equation*}
u_i^*(t)=\int_{-\infty}^t e^{-\be (t-s)} (f_i(u^*_s))+\be
u_i^*(s))\, ds,\q t\in \mathbb{R},
\end{equation*}
we obtain
\begin{equation*}
0=\int_{-\infty}^{t^*} e^{-\be (t^*-s)} (f_i(u^*_s))+\be u_i^*(s))\,
ds,
\end{equation*}
where $f_i(u^*_s)+\be u_i^*(s)\ge 0$ for $s\le t^*$. Hence
$f_i(u^*_s)+\be u_i^*(s)= 0$ for $s\le t^*$, and in particular
$u_i^*$ satisfies the scalar ODE $y'=-\be y$ for $s\le t^*$. Thus
$u_i^*(s)\equiv 0$ for $s\le t^*$, which is not possible. \smal

(iii) We note that  $u_t^*$ belongs to $W(0)$ for $t\le 0$, thus
$u_t^*=O(e^{\gamma t})$ at $-\infty$ and $u^*(t)$ satisfies
$u'(t)=Lu_t+h(t),$ with $h(t)=f(u_t^*)-Lu_t^*=O(e^{2\gamma t})$ at
$-\infty$. From Theorem \ref{TA.1} (see Appendix), for each $\vare
>0$ we deduce that $u^*(t)=z(t)+O(e^{(2\gamma -\vare) t})$ at
$-\infty$, where $z(t)=ce^{\la_0 t}{\bf v}$ for some (positive)
$c\in \mathbb{R}$. Thus, $u^*(t)=O(e^{\la_0 t})$ at $-\infty$. \ter

\begin{rem}\label{rem2.1} In fact,  one could use  \cite[Lemma 4]{FT} and
its constructive proof to derive that there is a complete solution
$u^*(t)$ of (\ref{e1.2}), with $u^*(-\infty )=0, u^*(\infty)=K$, and
$u^*(t)> 0$ for $t\le 0$. This proves assertion (i) of Theorem
\ref{T2.1}. In order to prove that  $u^*(t)=O(e^{\la_0 t})$ at
$-\infty$ it is however more convenient to explicitly construct
$u^*(t)$ as a perturbation of the eigenfunction $e^{\la_0 t}{\bf v}$
as above. This asymptotic result will be crucial to prove the
existence of {\it positive} travelling waves for (\ref{e1.1}), if
the wave speed is sufficiently high. On the other hand,   if we
assume that the interior of the positive cone $ {\cal C}_+$ is
positively invariant for the flow of (\ref{e1.2}), as an alternative
to hypothesis (H2)(ii), then the positiveness of $u^*(t)$ on
$\mathbb{R}$ follows immediately from the fact that $u^*(t)$ is
positive in the vicinity of $-\infty$.
\end{rem}

\section{Existence of travelling waves and their asymptotic
decay at $-\infty$} Throughout this section, for simplicity we
assume that (H1)-(H4) are fulfilled, but in fact some  of the
hypotheses can be weakened (cf. Remark \ref{R3.2}). We shall prove
the existence of  travelling waves for (\ref{e1.1}) which will be
obtained as perturbations of  the heteroclinic solution $u^*(t)$ of
(\ref{e1.2}). The asymptotic behaviour at $-\infty$ of $u^*(t)$
given in Theorem \ref{T2.1}(iii) will be important to study the
asymptotic decay of such waves  at $-\infty$; however its
 positiveness is irrelevant here, and will be only used for the analysis in Section 4.

For   a unit vector $w\in \mathbb{R}^p$,  we look for wave solutions
of (\ref{e1.1}) with direction $w$ and speed $c>0$, connecting the
equilibria 0 to $K$, i.e., solutions of the form $u(t,x)=\phi
(ct+w\cdot x)$ with $\phi(-\infty)=0,\phi(\infty)=K$.

The equation for the travelling wave profile $\phi$ is given by
\begin{equation}\label{e3.1}
\phi''(t)-c\phi'(t)+f_c(\phi_t)=0,\q t\in \mathbb{R},
\end{equation}
where $f_c(\phi)=f(\phi (c\cdot))$, with $\phi$ subject to  the
conditions
\begin{equation*}
\phi(-\infty)=0, \q \phi(\infty)=K.
\end{equation*}
With $\vare=1/c$, (\ref{e3.1}) is equivalent to
\begin{equation}\label{e3.2}
\vare^2\phi''(t)-\phi'(t)+f(\phi_t)=0.
\end{equation}
We also consider Eq. (\ref{e3.2}) with $\vare =0$, in which case it
reduces to Eq. (\ref{e1.2}).

Let $C_b(\mathbb{R},\mathbb{R}^N)$ be the space of all continuous
and bounded functions from $\mathbb{R}$ to $\mathbb{R}^N$, with the
supremum norm $\|y\|_{\infty}=\sup _{s\in \mathbb{R}}|y(s)|$. As a
particular case of the framework in \cite{FHW}, we have the
following result:
\begin{thm}\label{thm3.1}\cite{FHW} Let $f$ have the form $f(\phi)=F(\phi
(0),g(N\phi)),\, \phi\in {\cal C},$ for some  bounded linear
operator $N:{\cal C}\to \mathbb{R}^N$ and $g:\mathbb{R}^N\to
\mathbb{R}^N, \, F:\mathbb{R}^{2N}\to \mathbb{R}^N$ $C^2$-smooth
functions. Suppose also that:\vskip 0cm (i) $f(0)=f(K)=0$ for some
$K\in\mathbb{R}^N$,\vskip 0cm (ii) for Eq. (\ref{e1.2}), the
equilibrium $u=0$ is hyperbolic and unstable, and the equilibrium
$u=K$ is locally asymptotically stable. \vskip 0cm Then, if there is
a heteroclinic solution $u^*(t)$ for (\ref{e1.2}) connecting 0 to
$K$, for each unit $w\in \mathbb{R}^p$ there are a neighbourhood
${\cal V}$ of $u^*(t)$ in $C_b(\mathbb{R},\mathbb{R}^N)$ and  a
constant $c^*>0$, such that for $c>c^*$ the set of travelling waves
for (\ref{e1.1}) in ${\cal V}$, with direction $w$ and wave speed
$c$, constitutes a $C^1$-manifold of dimension $m$, where $m$ is the
dimension of the unstable space for $\dot u(t)=Df(0)u_t$. \end{thm}
In this section, the idea is to retrace some arguments in \cite{FHW}
for the proof of Theorem  \ref{thm3.1} adapted to the  case of
(\ref{e1.1}), but in appropriate Banach spaces, which will allow us
to deduce not only the existence of travelling wave solutions for
(\ref{e1.1}), but also their asymptotic behaviour  at $-\infty$.
This behaviour  will be used in Section 4, to prove the existence of
{\it positive} travelling waves.

\med In addition to  $C_b:=C_b(\mathbb{R},\mathbb{R}^N)$, we
introduce the following Banach spaces:

\smal

$C_b^1:=C_b^1(\mathbb{R},\mathbb{R}^N)=\{ y\in C_b:y'\in C_b\}$
with the norm $\|y\|_1=\|y\|_{\infty}+\|y'\|_\infty$;

$C_0=\{ y\in C_b:\lim _{s\to\pm \infty}y(s)=0\}$ is considered as a
subspace of $C_b$;

$C_0^1=\{ y\in C_b^1: y,y'\in C_0\}$ is considered as a subspace of
$C_b^1$;

$C_\mu=\{ y\in C_b:\sup_{s\le 0}e^{-\mu s}|y(s)|<\infty\}$ (for
$\mu>0$) with the norm
\begin{equation*}
\|y\|_\mu =\max \{\|y\|_\infty,\|y\|_{\mu}^- \}\q {\rm where} \q
 \|y\|_\mu^-=\sup_{s\le 0}e^{-\mu s}|y(s)|;
\end{equation*}

  $C_\mu ^1=\{ y\in C_b^1 : y,y'\in C_\mu\}$, with the norm $\|y\|_{1,\mu}=\|y\|_{\mu}+\|y'\|_\mu$;

  $C_{\mu,0}=C_{\mu}\cap C_0$ is considered as a subspace of $C_\mu$.

By the change of variables $\phi(t)=w(t)+u^*(t)$, (\ref{e3.2})
becomes
\begin{equation}\label{e3.3}
\vare^2w''(t)-w'(t)-w(t)=-w(t)-Df(u_t^*)w_t-G(\vare,t,w),\q t\in
\mathbb{R},
\end{equation}
where
\begin{equation}\label{e3.4}
G(\vare,t,w)=f(w_t+u_t^*)-f(u_t^*)-Df(u_t^*)w_t+\vare^2 {u^*}''(t),
\end{equation}
subject to the conditions $w(-\infty)=w(\infty)=0.$ The  roots of
the characteristic equation associated with
$\vare^2w''(t)-w'(t)-w(t)=0$ are
\begin{equation*}
\al(\vare)={{1-\sqrt{1+4\vare^2}}\over {2\vare^2}},\q
\be(\vare)={{1+\sqrt{1+4\vare^2}}\over {2\vare^2}},
\end{equation*}
and  satisfy $\al(\vare)\to -1^+,\be(\vare)\to \infty$ as $\vare \to
0^+$. In the case of different diffusion coefficients
$d_i>0,i=1,\dots,N,$ instead of $\al(\vare),\be(\vare)$ one has to
consider $\al_i(\vare),\be_i(\vare)$, the solutions of the
characteristic equations $d_i\vare^2z^2-z-1=0,i=1,\dots,N,$ but the
arguments are similar (cf. \cite{FHW}).

A bounded function $w:\mathbb{R}\to\mathbb{R}^N$ is a solution of
(\ref{e3.3}) if and only if
\begin{equation}\label{e3.5}
Jw(t)=H(\vare,w)(t),\q t\in\mathbb{R},
\end{equation}
where
\begin{equation*}
Jw(t)=w(t)-\int_{-\infty}^te^{-(t-s)} [w(s)+Df(u_s^*)w_s]\, ds
\end{equation*}
and
$$\displaylines{
H(\vare,w)(t)=\int_{-\infty}^t\left [{{e^{\al(\vare)(t-s)}}\over
{\sqrt{1+4\vare^2}}}-e^{-(t-s)}\right ](w(s)+Df(u_s^*)w_s)\, ds+\cr
{1\over {\sqrt{1+4\vare^2}}}\left [ \int_{-\infty}^t
e^{\al(\vare)(t-s)} G(\vare,s,w)\, ds+ \int^{+\infty}_t
e^{\be(\vare)(t-s)} [w(s)+Df(u_s^*)w_s+G(\vare,s,w)]\,
ds\right].\cr}$$

\smal

Our purpose is to apply a contraction principle argument in order to
obtain a solution of Eq. (\ref{e3.5}), for $\vare>0$ small and $w$
close to 0, in adequate spaces $C_\mu$. We first analyse the
linearity $J$, by introducing some auxiliary equations and
operators.

Define
\begin{equation*}
(Ty)(t)=y'(t)-Df(u_t^*)y_t,\q y\in C_b^1,t\in \mathbb{R}.
\end{equation*}
We easily see that  $J:C_0\to C_0,$ $T:C_b^1\to C_b$ are linear
bounded operators and $w\mapsto H( w,\vare)$ maps $C_0$ in $C_0$,
for $\vare>0$ (cf. \cite{FHW}). For $\mu>0$, we also define
\begin{equation*}
T_\mu:=T|_{C_\mu ^1}:C_\mu ^1\to C_\mu.
\end{equation*}
\begin{lem}\label{L3.1} Let $\mu>0$. Then $T_\mu$ and $J|_{C_{\mu
,0}}:{C_{\mu ,0}}\to {C_{\mu ,0}}$ are bounded  operators. \end{lem}
{\it Proof}. Since the map $t\mapsto \|Df(u_t^*)\|$ is continuous
for $t\in \R$ and $\|Df(u_t^*)\|\to \|Df(0)\|$ as $t\to -\infty$,
$\|Df(u_t^*)\|\to \|Df(K)\|$ as $t\to \infty$, then
$M:=\sup_{t\in\R} \|Df(u_t^*)\| <\infty.$
 It follows that
$|(Ty)(t)|\le \max (1,M)\|y\|_{1,\mu} $ for  $t\ge 0$ and $e^{-\mu
t}|(Ty)(t)|\le \max (1,M)\|y\|_{1,\mu} $ for $t\le 0, y\in C_\mu ^1
$.

For $y\in C_\mu$, we now have $|(Jy)(t)|\le (2+M)\|y\|_\mu$ for
$t\ge 0$ and $e^{-\mu t}|(Jy)(t)| \le [1+(1+M)(\mu
+1)^{-1}]\|y\|_\mu$, hence $J(C_{\mu ,0})\subset C_{\mu ,0}$.\ter

Consider the linear variational equation  around the heteroclinic
solution $u^*(t)$,
\begin{equation}\label{e3.6}
y'(t)=Df(u_t^*)y_t.
\end{equation}
Define the operators $L(t):=Df(u_t^*)$ in  (\ref{e3.6}), with
$L(-\infty)=Df(0)$ and $L(\infty)=Df(K)$. Hence, Eq. (\ref{e3.6}) is
asymptotically autonomous, with  limiting equations (\ref{e2.1}) and
$y'(t)=Df(K) y_t$, respectively at $-\infty$ and $\infty$.

\begin{lem}\label{L3.2} Consider $\mu\in (0,\la_0)$ such that there are no
characteristic roots $\la$ of (\ref{e2.1}) with $\Re\, \la =\mu$.
For $T_\mu:C_\mu ^1\to C_\mu$ defined as above,
\begin{equation*}
Im (T_\mu)=C_\mu,\q \dim\, Ker\, (T_\mu)=r_\mu,
\end{equation*}
where $r_\mu=\# \{ \la\in \mathbb{C} :\det \Delta_0(\la)=0, \Re\,
\la
>\mu\}.$ In particular, $r_\mu=1$ for $\mu$ close to $\la_0$.
Moreover, $Ker\, (T_\mu)\subset C_{\mu,0}$.
\end{lem}
{\it Proof}.
Clearly, equation $y'(t)=Df(K) y_t$ is asymptotically stable,
 and  the autonomous equation (\ref{e2.1})
admits  a ``shifted exponential dichotomy" in $\mathbb{R}$ with the
splitting made   at $\mu$ and exponents $\mu-\delta, \mu+\delta$,
for $\delta>0$ small. See Hale and Lin \cite{HL} for definitions,
and note that $C_\mu=C^0(\mu, 0)$ in the notation in \cite{HL}. From
\cite[Lemma 4.3]{HL}, there is $T>0$ such that  (\ref{e3.6}) has a
shifted exponential dichotomy on $(-\infty, -T]$ and $[T,\infty)$.
We now apply Lemma 4.6 of \cite{HL} to (\ref{e3.6}). It follows that
$T_\mu$ is a Fredholm operator, with index $Ind(T_\mu)$ given by
\begin{equation*}
Ind(T_\mu)=\dim Im(P_u^-(-t))-\dim Im (P_u^+(t)),\q t\ge T,
\end{equation*}
where $P_u^-(-t), P_s^-(-t)$ and $P_u^+(t), P_s^+(t)\, (t\ge T)$ are
the projections associated with the (shifted) exponential
dichotomies for $y'(t)=Df(0) y_t$ and $y'(t)=Df(K) y_t$,
respectively. From \cite[Lemma 4.3]{HL}, we also have that
$P_u^-(-t)\to P_u^-, P_u^+(t)\to P^+_u$ as $t\to\infty$, where
$P_u^-$ is the canonical projection from ${\cal C}$ onto the
$\mu$-unstable space $E_\mu^-$ for   $y'(t)=Df(0) y_t$, and $P_u^+$
is the canonical projection from ${\cal C}$ onto the unstable space
$E_u^+$ for for $y'(t)=Df(K) y_t$. We have $E_u^+=\{ 0\}$ and $\dim
E_\mu^-=r_\mu$, where $r_\mu$ is the number of characteristic values
for (\ref{e2.1}) (counting multiplicities) with real parts greater
than $\mu$. Hence $Ind (T_\mu)=r_\mu$.  On the other hand, the index
of $T_\mu$ is defined by $Ind (T_\mu)=\dim Ker(T_\mu)-{\rm codim}\,
Im(T_\mu)$. Again by \cite[Lemma 4.6]{HL} we have $\dim
Ker(T_\mu)=\dim E_\mu^-=r_\mu$, yielding that $Im (T_\mu)=C_\mu$.

For  $y\in Ker\, (T_\mu)$, from the definition of shifted
exponential dichotomy we have $\lim_{t\to\infty} y(t)=0.$ Thus,
$Ker\, (T_\mu)\subset C_{\mu,0} $. \ter

Similarly to what was done for $T$, we now restrict the domain and range of the operator $J$.
With $D:=d/dt +id$, consider the commutative diagram
\begin{equation*}
\begin{array}[c]{ccccc}
C^1_\mu  &&\stackrel{T_\mu}\longrightarrow&&C_\mu\\
& {_{J}}{\searrow}&&{\nearrow}_D\\
&&C^1_\mu
\end{array}
\end{equation*}
It is easy to check that this diagram is well defined, and that $D$
is one-to-one and surjective. Since $T_\mu = D\circ  J$ is
surjective, we may conclude that $J$ is also surjective. Moreover,
\begin{lem}\label{L3.3} Consider $\mu\in (0,\la_0)$ such that there are no
characteristic roots $\la$ of (\ref{e2.1}) with $\Re\, \la =\mu$.
Then, for the operator $J|_{C_{\mu ,0}}:{C_{\mu ,0}}\to {C_{\mu
,0}}$  we have $Ker\, (J|_{C_{\mu ,0}})=Ker\, (T_\mu)$ and $Im\,
(J|_{C_{\mu ,0}})=C_{\mu ,0}$.
\end{lem}

{\it Proof}. Recall that $Ker\, (T_\mu)\subset C_{\mu,0}$.
 Clearly, for $w\in C_{\mu ,0}$ we have $Jw=0$ if and only if $w'(t)=Df(u_t^*)w_t$, and   then
 $w'\in C_\mu$. We therefore deduce  that $(Ker\, J)\cap C_\mu =(Ker\, J)\cap C_\mu^1$, and $Ker\, (J|_{C_{\mu ,0}})=Ker\, (T_\mu)$.

We now prove that $Im\, (J|_{C_{\mu ,0}})=C_{\mu ,0}$. Indeed, for
$y \in C_{\mu,0}$ we have that  $\xi:= y-Jy \in C^1_{\mu}$ and $D\xi
(t)=y(t)+Df(u_t^*)y_t$, hence $D\xi \in C_{\mu,0}$. Equation $Jw =
y$ is equivalent to $ J(w-y) =  \xi ,$ and therefore it possesses a
solution $\chi= w-y \in C^1_\mu$. After applying $D$ to both sides
of the latter equation, we get $ T_\mu \chi = D\xi \in C_{\mu,0}$.
Since the $\omega$-limit operator $T_\mu(\infty)$ is hyperbolic, we
may invoke Lemma 3.3 from \cite{FHW} to conclude that
$\chi(\infty)=0$. Thus $w(\infty)= 0$, and $J: C_{\mu,0} \to
C_{\mu,0}$ is surjective. \ter

 We now focus our attention on the non-linearity $H$ of Eq. (\ref{e3.5}). Proceeding as in \cite{FHW},
one sees  that $H(w,\vare)\in C_0$ for each $\vare > 0$ and $w\in
C_0$. We  want however to  consider the maps $H(\cdot ,\vare)$
restricted to some neighbourhood of zero in $C_{\mu,0}$, for
$\vare>0$ and $\mu >0$. We start with an auxiliary lemma:

\begin{lem}\label{L3.4} Let  $X, Y$ be normed spaces and $C\subset O \subset X$.
Suppose that the set $C$ is compact, $O$ is open and $F:O \to Y$ is
a continuous map. Then for every $\sigma >0$ there exists $\delta
>0$ such that
\begin{equation*}
|F(x+z)-F(x)| \leq \sigma, \quad x \in C, \  |z| \leq \delta.
\end{equation*}
\end{lem}
{\it Proof}. By the continuity of $F$, for each $x\in C$ there is
$\delta(x)>0$ such that if $|z| \leq 2\delta(x)$ then $x+z \in O$
and $ |F(x+z)-F(x)| \leq \sigma/2$. Since $C \subset \cup_{x \in
K}B_{\delta(x)}(x)$ is compact, there is a finite subcover
$\{B_{\delta(x_j)}(x_j)\}_{j=1}^m$ of $C$. For each $x \in
B_{\delta(x_j)}(x_j)\cap C$ and $|z|\leq \delta := \min
\{\delta(x_j)\}$, we have $ |F(x+z)-F(x)| \leq |F(x+z)-F(x_j)|+
|F(x)-F(x_j)| \leq \sigma,$ which proves the lemma. \ter

\begin{lem}\label{L3.5} Assume (H1)--(H4) and  consider $\mu\in (0,\la_0)$.
  For any $\delta>0$, there are  $\vare ^*>0$ (independent of $\mu$)  and $\sigma >0$  such that $H(w,\vare)\in C_{\mu ,0}$
   for any $\vare >0$ and  $w\in C_{\mu ,0}\cap B_\sigma^\mu(0)$,    and
   \begin{eqnarray}\label{e3.7}
  \nonumber \fl  \|H(w,\vare)\|_\mu \le  \delta (\|w\|_\mu +1), && \\
    \fl  \|H(w,\vare)-H(v,\vare)\|_\mu \le \delta \|w-v\|_\mu, &&\q {\rm for}\q w,v\in C_{\mu,0} \cap
B_\sigma^\mu(0), \vare\in (0,\vare^*)
   \end{eqnarray}
 where $B_\sigma^\mu(0)$ is the $\sigma$-neighbourhood of 0 in $C_\mu$.
\end{lem}
{\it Proof}.
 We write $H=H_1+H_2+H_3$, where
 $$\eqalign{
H_1(\vare,w)(t)&=\int_{-\infty}^t\left [{{e^{\al(\vare)(t-s)}}\over
{\sqrt{1+4\vare^2}}}-e^{-(t-s)}\right ](w(s)+Df(u_s^*)w_s)\, ds,\cr
H_2(\vare,w)(t)&={1\over {\sqrt{1+4\vare^2}}} \int_{-\infty}^t
e^{\al(\vare)(t-s)} G(\vare,s,w)\, ds,\cr H_3(\vare,w)(t)&={1\over
{\sqrt{1+4\vare^2}}}\int^{+\infty}_t e^{\be(\vare)(t-s)}
[w(s)+Df(u_s^*)w_s+G(\vare,s,w)]\, ds,\cr }$$ and $G$ is given by
(\ref{e3.4}). Let $M=\sup_{t\in\R}\|Df(u_t^*)\|$ as before. For
$t\in\mathbb{R},\ \vare> 0$ and $\mu \ge 0$, we have
$$\eqalign{
&\left |\int_{-\infty}^t\left [{{e^{\al(\vare)(t-s)}}\over
{\sqrt{1+4\vare^2}}}-e^{-(t-s)}\right ]e^{\mu s}\, ds\right |\cr
&\le {1\over {\sqrt{1+4\vare^2}}}\int_{-\infty}^t \left |
e^{\al(\vare)(t-s)} (1-\sqrt{1+4\vare^2})+\sqrt{1+4\vare^2}\Big
(e^{\al(\vare)(t-s)}-e^{-(t-s)}\Big )\right | e^{\mu s} \, ds \cr &=
{1\over {\sqrt{1+4\vare^2}}}\left [
(\sqrt{1+4\vare^2})-1)\int_{-\infty}^t  e^{\al(\vare)(t-s)} e^{\mu
s} \, ds  + \sqrt{1+4\vare^2} \int_{-\infty}^t
(e^{\al(\vare)(t-s)}-e^{-(t-s)}\Big ) e^{\mu s} \, ds \right ] \cr
&=  {1\over {\sqrt{1+4\vare^2}}}\left [
 {{2\sqrt{1+4\vare^2}-1}\over {\mu -\al(\vare)}}- { {\sqrt{1+4\vare^2}}\over {\mu +1}}\right] e^{\mu
 t}\cr}
 $$
 \begin{equation}\label{e3.8}
 \fl = {1\over {\sqrt{1+4\vare^2}}} \left [ {{\sqrt{1+4\vare^2}-1}\over
{\mu -\al(\vare)}}+  {{(1+\al(\vare))\sqrt{1+4\vare^2}}\over
{(\mu-\al(\vare))(\mu+1)}} \right ] e^{\mu t}
 \le C_1(\vare)e^{\mu t},
 \end{equation}
where
\begin{equation*}
C_1(\vare)= -{1\over {\al(\vare)}} \left ( 1- {1\over
{\sqrt{1+4\vare^2}}}+1+\al(\vare)\right ) \to 0\q {\rm as}\q \vare
\to 0^+.
\end{equation*}
From (3.8), we obtain
\begin{equation}\label{e3.9}
\fl \|H_1(\vare,w)-H_1(\vare, v)\|_\mu\le C_1(\vare)
(1+M)\|w-v\|_\mu,\q w,v\in C_\mu,\vare> 0.
\end{equation}
Since $H_1(\vare, 0)=0$, in particular $H_1(\vare,w)\in C_\mu$ for
$w\in C_\mu$ and $\vare >0$.

For $0\le \mu<\be (\vare)$ and $t\in\mathbb{R}$, we now have
\begin{equation}\label{e3.10}\fl
\int_{-\infty}^t e^{\al(\vare)(t-s)}e^{\mu s}\, ds={{e^{\mu t}}\over
{\mu-\al(\vare)}},\ \int^{+\infty}_t e^{\be(\vare)(t-s)} e^{\mu s}
\, ds={{e^{\mu t}}\over {\be (\vare)-\mu}}.
\end{equation}
Consider e.g. $\mathbb{R}^N$ equipped with the maximum norm. For
$t\in\mathbb{R},$ $\vare > 0,$ $ w,v\in C_{\mu,0},$ $ i=1,\dots,N$,
we have
\begin{eqnarray}\label{e3.11}
\nonumber \fl &  |G_i(\vare, t,w)|&\le
\vare^2|{u_i^*}''(t)|+|f_i(w_t+u_t^*)-f_i(u_t^*)-Df_i(u_t^*)w_t| \\
\fl  & & \le
\vare^2|{u_i^*}''(t)|+\|Df_i(u^*_t+\xi_{i,t}w_t)-Df_i(u_t^*)\|\|w_t\|_\infty
\end{eqnarray}
and
\begin{equation}\label{e3.12}\fl
|G_i(\vare, t,w)-G_i(\vare,t,v)|\le
\|Df_i(v_t+u_t^*+\th_{i,t}(w_t-v_t))-Df_i(u_t^*)\|\|w_t-v_t\|_\infty,
\end{equation}
for some $\xi_{i,t},\th_{i,t} \in (0,1)$ for $t\in\mathbb{R}$.

Note  that ${u^*}''\in C_\mu$ for $0<\mu \le \la_0$. In fact,
$u^*\in C_\mu$  from Theorem \ref{T2.1}, hence Eq. (\ref{e1.2}) and
the smoothness of $f$ lead to $|{u^*}'(t)|\le M_0\|u^*_t\|_\infty$,
from which we derive $\|{u^*}'\|_\mu \le M_0 \|u^*\|_\mu,$ for some
$M_0>0$. By differentiating, we obtain ${u^*}''(t)=Df(u_t^*)
({u^*}')_t$,  thus $\|{u^*}''\|_\mu \le M \|{u^*}'\|_\mu.$

In order to simplify the notation, for each $\mu ,\sigma >0$ write
$C_{\mu,0}\cap B_{\sigma}(0)$ to denote the $\sigma$-neighbourhood
of 0 in $C_{\mu ,0}$. Since $u^*$ is  uniformly bounded on
$\mathbb{R}$ and $f$ transforms bounded sets of ${\cal C}$ into
bounded set of $\mathbb{R}^N$, then ${u^*} '$ is uniformly bounded
on $\mathbb{R}$ and $u^*$ uniformly continuous on $\mathbb{R}$.
Thus, ${\cal K}= \overline{\{{u_t^*}, \ t \in \mathbb{R}\}} \subset
{\cal C}$ is compact. The continuity of $Df_i: {\cal C} \to
\cal{L}(C,C)$ and Lemma \ref{L3.4} imply that,  for each $\delta>0$
fixed, there is $\sigma =\sigma(\de,\mu)>0$ such that
$\|Df_i(v_t+u_t^*+\th_{i,t}(w_t-v_t))-Df_i(u_t^*)\|<\delta$ for
$w,v\in C_{\mu,0}\cap B_\sigma(0), \ t \in \R$. From (\ref{e3.11}),
(\ref{e3.12}), we get
$$\displaylines{
|G(\vare, t,w)|\le \vare^2|{u^*}''(t)|+\delta \|w_t\|_\infty, \q
w\in C_{\mu,0}\cap B_{\sigma}(0),\cr |G(\vare, t,w)-G(\vare,t,v)|\le
\delta  \|w_t-v_t\|_\infty,\q w,v\in C_{\mu,0}\cap B_\sigma(0).\cr}
$$
 From these estimates and (\ref{e3.10}), we conclude that
$H_2(\vare,w), H_3(\vare,w)\in C_\mu$ for all $w\in C_\mu\cap
B_{\sigma}(0)$ and $\vare >0,\mu \in (0,\la_0)$, with
\begin{eqnarray}\label{e3.13}
\nonumber \fl & & \|H_2(\vare,w)\|_\mu \le {1\over
{(\mu-\al(\vare))\sqrt{1+4\vare^2}}}(\vare^2\|{u^*}''\|_\mu+\delta
\|w\|_\mu), \\
\fl  & & \|H_3(\vare,w)\|_\mu \le {1\over
{(\be(\vare)-\mu)\sqrt{1+4\vare^2}}} \Big
[\vare^2\|{u^*}''\|_\mu+(1+M+\delta )\|w\|_\mu \Big].
\end{eqnarray}
Furthermore, for $\|w-v\|_\mu, w,v\in C_{\mu,0}\cap
B_\sigma(0),\vare>0,$ we get
\begin{eqnarray}\label{e3.14}
\nonumber \fl & & \|H_2(\vare,w)-H_2(\vare, v)\|_\mu\le{\delta \over
{(\mu-\al(\vare))\sqrt{1+4\vare^2}}} \|w-v\|_\mu, \\
\fl  & & \|H_3(\vare,w)-H_3(\vare, v)\|_\mu\le {{1+M+\delta }\over
{(\be(\vare)-\mu)\sqrt{1+4\vare^2}}}.
\end{eqnarray}
 On the other hand,
$${1\over {\sqrt{1+4\vare^2}}}\left ({1\over {\be(\vare)-\mu}}+{1\over {\mu -\al(\vare)}} \right ) ={1\over {1+\mu-\vare^2\mu^2}}< 1$$
if $\vare^2\mu <1$. From (\ref{e3.9}), (\ref{e3.13}) and
(\ref{e3.14}), for $\vare>0$ small enough and $\mu\in (0,\la_0)$ we
obtain
\begin{equation}\label{e3.15}
\|H(w,\vare)\|_\mu\le C(\vare)\|w\|_\mu +D(\vare),\q w\in
C_{\mu,0}\cap B_{\sigma}(0)
\end{equation}
and
\begin{equation}\label{e3.16}
\|H(w,\vare)-H(v,\vare)\|_\mu\le C(\vare) \|w-v\|_\mu,\q w,v\in
C_{\mu,0} \cap B_\sigma(0),
\end{equation}
where $C(\vare), D(\vare)$ do not depend on $\mu$ and are given by
$$C(\vare)=C_1(\vare)(1+M)+\delta+
{{1+M }\over {(\be(\vare)-\la_0)\sqrt{1+4\vare^2}}},\q D(\vare)=
\vare^2\|{u^*}''\|_{\la_0} .$$ Since  $C_1(\vare)\to 0,
\be(\vare)\to\infty$ as $\vare\to 0^+$, by replacing $\delta$ by
$\delta /2$ in  (\ref{e3.15}), (\ref{e3.16}), we obtain (\ref{e3.7})
for $\vare>0$ sufficiently small. \ter

We now  return to Eq. (\ref{e3.5}).  Let $0< \mu <\la_0$. For $\vare
>0$ small, we look for a solution $w\in C_{\mu ,0}$ of (\ref{e3.5}). For
the case $\mu=0$, where the space $C_{0,0}$ denotes $C_0$,  this
question was addressed in \cite{FHW}. Our purpose is to solve this
problem for $\mu\in (0,\la_0)$.

We first apply a Liapunov-Schmidt reduction. From Lemmas \ref{L3.2}
and \ref{L3.3}, $X_\mu:=Ker\, (J|_{C_{\mu ,0}})$ is finite
dimensional, hence there is a complementary subspace $Y_\mu$ in
$C_{\mu ,0}$,
\begin{equation*}
C_{\mu ,0}=X_\mu \oplus Y_\mu.
\end{equation*}
For $w\in C_{\mu ,0}$, write $w=\xi+\phi$ with $\xi\in X_\mu,\phi\in
Y_\mu$. Define $S_\mu:=J|_{Y_\mu}$. Since $S_\mu:Y_\mu\to C_{\mu,0}$
is bounded and bijective, then $S_\mu^{-1}$ is bounded. In the space
$C_{\mu,0}$, Eq (\ref{e3.5}) is equivalent to $\phi
=S_\mu^{-1}H(\vare, \xi+\phi)$, therefore we look for fixed points
$\phi\in Y_\mu$ of the map
\begin{equation}\label{e3.17}
{\cal F}_\mu(\vare,\xi,\phi)=S_\mu^{-1}H(\vare, \xi+\phi).
\end{equation}
For simplicity, in what follows we write $S, {\cal F}, B_\sigma(0)$
instead of $S_\mu,{\cal F}_\mu,B_\sigma^\mu (0)$, respectively, when
there is no risk of misunderstanding.
\begin{rem}\label{R3.1}
For $0<\mu_1<\mu_2 <\la_0$ with $\mu_1,\mu_2\notin \Re\, \sigma
(A)$, where $ \sigma (A)$ is the set of solutions of (\ref{e2.2}),
it is clear that $C_{\mu_2}\subset C_{\mu_1}$ with $\|y\|_{\mu_1}\le
\|y\|_{\mu_2}$, and $X_{\mu_2}\subset X_{\mu_1}$. Together with
Lemmas \ref{L3.2} and \ref{L3.3}, this implies that for each
interval $I:=[\mu_1,\mu_2] \subset (0,\la_0)\setminus \Re\, \sigma
(A)$, we have $X_{\mu_2}= X_{\mu_1}$. We now show that the
complementary subspaces $Y_\mu$ can be chosen so that
$Y_{\mu_2}\subset Y_{\mu}\subset Y_{\mu_1}$ for $\mu\in I$. In fact,
let $X_\mu=span\, \{ y_1,\dots, y_r\}$, where $y_1,\dots, y_r\in
C_{\mu,0}$ and $r=r_\mu$ for $\mu\in I$. From the Hahn-Banach
theorem, let $h_i\in (C_{\mu_1,0})'$ be such that $h_i(y_i)=1,
h_i(y_j)=0$ for $j\ne i, i,j=1,\dots, r$. Define the natural
injections $i(\mu,\mu_1): C_{\mu,0}\to C_{\mu_1,0}$, which are
continuous, and the subspaces $Y_\mu = \{ y\in C_{\mu,0}: h_i\circ
i(\mu,\mu_1)(y)=0, i=1,\dots, r\}$. Hence $Y_\mu$ is a closed
subspace of $C_{\mu,0}$, and  for $y\in C_{\mu,0}$ we have
$\sum_{i=1}^r h_i(y)y_i\in X_\mu, y-\sum_{i=1}^r h_i(y)y_i\in
Y_\mu$, from which  the decompositions $C_{\mu ,0}=X_\mu \oplus
Y_\mu$ follow, with $Y_{\mu_2}\subset Y_{\mu}\subset Y_{\mu_1}$ for
$\mu\in I$.
\end{rem}

\begin{thm} \label{Thm3.2} Assume (H1)-(H4),
and denote by $\sigma(A)$ the set of characteristic values for
(\ref{e2.1}). Fix an interval $I:=[\mu_1\mu_2]\subset (0,\la_0)
\setminus \Re\, \sigma (A)$, and denote $r=r_\mu$ for all $\mu\in
I$. Then, there exist $\vare^*>0$ and $\sigma >0$, such that for
$0<\vare\le \vare^*$, the following  holds: for each unit vector
$w\in \mathbb{R}^p$ and all $\mu\in I$, in a neighbourhood
$B_\sigma^\mu(0)$ of $u^*(t)$ in $C_\mu$, the set of all travelling
wave solutions $u(t,x)=\psi (ct+w\cdot x)$ of (\ref{e1.1}) with
speed $c=1/\vare$ and connecting 0 to $K$ forms a $r$-dimensional
manifold (which does not depend on $\mu$), with the profile
 $\psi\in {\cal M}_{I,\vare}$, where
 $${\cal M}_{I,\vare}=\{ \psi: \psi (t)=u^*(t)+\xi +\phi (\vare, \xi),\ {\rm for}\ \xi\in X_\mu\cap B_\sigma ^\mu(0)\} ,$$
where  $\phi (\vare, \xi)=\phi_\mu (\vare, \xi)$ is the fixed point
of ${\cal F}_\mu(\vare,\xi,\cdot)$ in $Y_\mu\cap B_\sigma^\mu (0)$,
and is continuous on $(\vare,\xi)$.
\end{thm}

{\it Proof}.   In the sequel, we shall use the simplified notation
$S, {\cal F}, B_\sigma(0)$,   for $S_\mu, {\cal F}_\mu,B_\sigma^\mu
(0)$, respectively. Fix $\mu \in I$ and $k\in (0,1)$. From Lemma
\ref{L3.5}  (cf. (\ref{e3.15}) and (\ref{e3.16})), for $\delta>0$
small there are $\sigma=\sigma(\de,\mu)>0$ and
$\vare^*=\vare^*(\de)>0$ such that for $0<\vare<\vare^*$, $\xi \in
X_\mu\cap \overline{B_\sigma (0)}$ and $\phi_1,\phi_2\in Y_\mu\cap
\overline{B_\sigma (0)}$ we have
\begin{equation}\label{e3.18} \fl
 \|S^{-1}H(\xi+\phi_1,\vare)\|\le  \|S^{-1}\| \big (C(\vare)\|\xi+\phi_1\|_\mu +D(\vare)\big )\le \delta \|S^{-1}\| (\|\xi+\phi_1\|_\mu+1)
\end{equation}
and
\begin{equation}\label{e3.19} \fl
 \|S^{-1}(H(\xi+\phi_1,\vare)-H(\xi+\phi_2,\vare))\|_\mu \le \|S^{-1}\| C(\vare)  \|\phi_1-\phi_2\|_\mu \le \delta\|S^{-1}\| \|\phi_1-\phi_2\|_\mu
\end{equation}
with $\delta (1+2\sigma)\|S^{-1}\|\le \sigma$ and $\delta \|S^{-1}\|
\le k$.   From (\ref{e3.18}) and (\ref{e3.19}), it follows that
${\cal F}:(0,\vare^*)\times (X_\mu\cap \overline{B_\sigma
(0)})\times (Y_\mu\cap \overline{B_\sigma (0)})\to Y_\mu\cap
\overline{B_\sigma (0)}$  is a uniform contraction map of $\phi \in
Y_\mu\cap \overline{B_\sigma (0)}$, hence for $(\vare, \xi)\in
(0,\vare^*)\times (X_\mu\cap \overline{B_\sigma (0)})$ there is a
unique solution $\phi (\vare, \xi)=\phi_\mu (\vare, \xi)\in Y_\mu$
of (\ref{e3.17}), with $\phi (\vare, \xi)$ continuous. Define the
$r$-dimensional manifold ${\cal M}_{\mu,\vare}=\{ \psi: \psi
(t)=u^*(t)+\xi +\phi (\vare, \xi),\ {\rm for}\ \xi\in X_\mu\cap
B_\sigma (0)\} $. Choose $\sigma=\sigma(\de,\mu_2)$, independent of
$\mu\in I$. From the uniqueness of the fixed point and Remark
\ref{R3.1}, it follows that $\phi_\mu (\vare, \xi)=\phi_{\mu_2}
(\vare, \xi)$ does not depend on $\mu\in I$, as well as  ${\cal
M}_{\mu, \vare}:={\cal M}_{I, \vare}$.\ter

We observe that if 0 is a hyperbolic equilibrium of (\ref{e1.2}) and
$f$ has the particular form $f(\phi)=F(\phi(0), g(N\phi))$, then
Theorem \ref{thm3.1} asserts that the result in Theorem \ref{Thm3.2}
is valid for $\mu=0$.
\begin{cor}\label{Cor3.1} Under the assumptions of Theorem \ref{Thm3.2} and
with the same notation, for $0<\mu<\la_0 $ such that the strip $\{
\la\in \mathbb{C}: \Re\, \la \in (\mu, \la_0)\}$ does not intersect
$\sigma (A)$, the manifold ${\cal M}_{\mu, \vare}$ is 1-dimensional.
\end{cor}
\begin{cor}\label{Cor3.2} Under the assumptions of Theorem
\ref{Thm3.2} and with the same notation, for  an interval
$I:=[\mu_1\mu_2]\subset (0,\la_0) \setminus \Re\, \sigma (A)$, there
are $\vare^*>0,\sigma
>0$ and $C>0$ such that  the travelling profiles $\psi (\vare,\xi)$
satisfy
\begin{equation*}
\|\psi (\vare,\xi)\|_\mu \le C, \q \|\psi '(\vare,\xi)\|_\mu \le C\q
{\rm for}\q 0<\vare<\vare^*, \xi \in X_\mu\cap \overline{B_\sigma
(0)},
\end{equation*}
where $C$ does not depend on  $\mu\in I$. In particular $|\psi
(\vare, \xi)(t)|\le Ce^{\mu t}$ for $t\le 0, 0<\vare<\vare^*,$ $ \xi
\in X_\mu\cap \overline{B_\sigma (0)}.$ \end{cor}

{\it Proof}. For all $\mu\in I$, the profiles  are given by  $\psi
(\vare,\xi)=u^*+\xi+ \phi (\vare,\xi)$, where $\phi
(\vare,\xi)=\phi_{\mu_2}  (\vare,\xi)$ is the fixed point of $ {\cal
F}_{\mu_2}(\vare,\xi,\cdot).$ Since $\|y\|_\mu \le \|y\|_{\mu_2}$
for $y\in C_{\mu_2,0}$, we only need to prove the result for
$\mu=\mu_2$. In what follows, we write $S_{\mu_2}^{-1}=S^{-1}$.

Fix $k\in (0,1)$, and consider $\vare_1>0$ such that $\|S^{-1}\|
C(\vare) \le k$ for $0<\vare<\vare_1$ where $C(\vare)$ is as in
(\ref{e3.16}). From (\ref{e3.19}), for $w_1,w_2\in C_{\mu_2 ,0}\cap
B_\sigma (0)$ we have
$$
\fl \|S^{-1}(H(w_1,\vare)-H(w_2,\vare))\|_{\mu_2} \le k
\|w_1-w_2\|_{\mu_2},
$$
and  the contraction principle yields
\begin{equation*}
\fl \|\phi (\vare,\xi)\|_{\mu_2} \le {1\over {1-k}} \|  {\cal
F}_{\mu_2}(\vare,\xi,0)\|_{\mu_2} \le {1\over {1-k}} \|S^{-1}\|
\|H(\vare,\xi)\|_{\mu_2},\q \xi \in X_{\mu}\cap \overline{B_\sigma
(0)}.
\end{equation*}
For $\vare,\sigma >0$ small enough, from (\ref{e3.15}) we get
$$\| H(\vare, \xi)(t)\|_{\mu_2} \le C(\vare)\|\xi\|_{\mu_2} +\vare^2 \| {u^*}''\|_{\mu_2},\q
 \xi \in X_{\mu}\cap \overline{B_\sigma (0)}.$$
We thus  obtain $\| \psi (\vare,\xi)\|_{\mu_2} \le
\|u^*\|_{\mu_2}+(1+ C(\vare))\sigma +\vare^2 \| {u^*}''\|_{\mu_2}\le
C_1$ for  $\vare$ small, where $C_1$ does not depend on $\vare,\xi$.

Now we want to prove a similar estimate for the derivates
$d\psi(\vare,\xi)/dt$. For simplicity, we only prove the result for
$\xi=0$.

Since $\psi (t):=\psi (\vare,0)(t)$ is a solution of (\ref{e3.2}),
then $\psi (t)$ is given by the integral formula
\begin{equation*}
 \fl \psi (t)={1\over {\sqrt{1+4\vare^2}}}\left ( \int_{-\infty}^t
e^{\al(\vare)(t-s)} [\psi (s)+f(\psi_s)]\, ds+ \int^{+\infty}_t
e^{\be(\vare)(t-s)}   [\psi (s)+f(\psi_s)]\, ds\right ),
\end{equation*}
from which we derive
\begin{equation*}\fl \psi'(t)={1\over {\sqrt{1+4\vare^2}}}\Big
(\al(\vare)\hspace{-2mm} \int\limits_{-\infty}^t
e^{\al(\vare)(t-s)}[\psi (s)+f(\psi_s)]ds
 -\be(\vare)\hspace{-2mm}\int\limits^{+\infty}_t e^{\be(\vare)(t-s)}[\psi (s)+f(\psi_s)]
ds\Big ).
\end{equation*}
Since  $f$ is bounded on bounded sets of ${\cal
C}=C([-\tau,0];\mathbb{R}^N)$, there is $\ell$ such that
$|f(\phi)|\le \ell $ for $\phi\in {\cal C}$ with $\|\phi\|_\infty
\le C_1$. Thus, $|f(\psi_s)|\le \ell $ for $s\in\mathbb{R}$, where
$\ell$ does not depend on $\mu,\vare$, and $\|\psi'\|_\infty\le
2(C_1+\ell)/\sqrt{1+4\vare^2}.$ From (\ref{e3.10}), the
$C^1$-smoothness of $f$ and $f(0)=0$, we easily deduce that there is
$C_2>0$ such that $\|\psi'\|_\mu \le C_2$. This completes the
proof.\ter

In fact a stronger result can be proven:

\begin{cor} \label{Cor3.3} Assume (H1)-(H4), take $\mu \in  I:=[\mu_1\mu_2]\subset (0,\la_0) \setminus \Re\, \sigma (A)$,
 and consider  the travelling wave profiles $\psi (\vare,\xi)=u^*+\xi+\phi (\vare, \xi)$
for $\vare\in (0,\vare^*),\xi\in X_\mu\cap B_\sigma (0)$, given in
Theorem \ref{Thm3.2}. For $\xi=0$ and $\mu\in I$, the profile $\psi
(\vare,0)$  satisfies
\begin{equation}\label{e3.20}
\psi (\vare,0)\to u^*\ {\rm in} \ C_\mu\q {\rm as} \q \vare\to 0^+.
\end{equation}
\end{cor}
 {\it Proof.}  Let $\vare^*,\sigma>0$ be as in the statement of
Theorem \ref{Thm3.2}, and recall that  $\psi (\vare,0)=
\psi_\mu(\vare,0)$ only depends on $\vare$. Next, we deduce some
estimates  as in Lemma \ref{L3.5}, so details are omitted.
   For $\vare=0$, define
    $$H(0,w)(t)=\int_{-\infty}^t e^{-(t-s)} [f(w_s+u_s^*)-f(u_s^*)-Df(u_s^*)w_s]\, ds.$$
     We write  $H(0,w)(t)=\int_{-\infty}^t e^{-(t-s)} G(0,t,w)\, ds:=H_2(0,w)$, where $G(0,t,w)$ is given by (\ref{e3.4}).
     After some computations, we observe  that the function $H$ restricted to $[0,\vare^*)\times (C_{\mu,0}\cap B_\sigma^\mu (0))$ satisfies
   $$\|H(\vare,w)-H(0,w)\|_\mu\le C_0(\vare)\|w\|_\mu+D_0(\vare),
   $$
   with   $C_0(\vare),D_0(\vare)$ independent of $\mu$, $C_0(\vare),D_0(\vare)\to 0$ as $\vare\to 0^+$.
   This means that the function $(\vare,w)\mapsto H(\vare,w)$ converges, uniformly on $w\in C_{\mu,0}\cap B_\sigma^\mu (0)$,   to $H(0, \cdot )$ in  $C_{\mu,0}$ as $\vare \to 0^+$.

Moreover, for $\vare =0$ and  $\xi=0$  the fixed point of (3.17) is
$\phi (0,0)=0$. Therefore, the application of the contraction
principle as in the proof of Theorem \ref{Thm3.2} leads to
(\ref{e3.20}).\ter
\begin{rem}\label{R3.2} As seen in Section 2, the existence of a positive
eigenvector ${\bf v}\in\mathbb{R}^N$ associated with the
characteristic root $\la_0$ of (\ref{e2.2}) was crucial to prove the
existence of a {\it positive} heteroclinic solution $u^*(t)$ of
(\ref{e1.2}), connecting the equilibria 0 to $K$. For all the
results in this section the positiveness of such heteroclinic
solution is irrelevant, and therefore it is not necessary to impose
the above requirement  in (H4) that ${\bf v}$ is positive. For the
same reason, in Section 3 assumption (H2)(ii) is not needed as well.
\end{rem}

\section{Positiveness of travelling waves}
Consider the characteristic equation for the linearization of
(\ref{e3.2}) at 0,
\begin{equation}\label{e4.1}
\det \Delta_\vare (z)=0,\q {\rm where}\q \Delta_\vare (z):=\vare^2
z^2I-zI+L(e^{z\cdot}I),
\end{equation}
where $L=Df(0)$.  Recall that for $\vare=0$ the characteristic
matrix-valued function $\Delta _0(z)$ was defined in (\ref{e2.2}).
Since $\la_0>0$ is a simple root of the characteristic equation
$\det \Delta_0 (z)=0$, from the implicit function theorem, for
$\vare >0$ small there is a simple real  root $\la(\vare)$ of
(\ref{e4.1}), with $\la(\vare)\to \la_0$ as $\vare\to 0^+$.

\begin{lem}\label{L4.1} For $\de >0$ sufficiently small and $\de_1
>0$,  there exists $\vare_0>0$ such that,   for $0<\vare<\vare_0$,
$\la(\vare)$ is the only root of the characteristic equation
(\ref{e4.1}) on the vertical strip $\la_0-\de \le Re\, z\le \la_0
+\delta_1$.
\end{lem}
{\it Proof}. Let  $\delta >0$ be such that $\la_0$ is the only root
of  $\det \Delta_0 (z)=0$ on the strip $S=\{ z:\la_0-\delta \le
\Re\, z\le \la_0+\delta_1\}$. If $z(\vare)\in S$ is a root of $\det
\Delta_\vare (z)=0$, then there is a unit vector
$w=w(\vare)\in\mathbb{R}^N$ such that
$(\vare^2z(\vare)^2-z(\vare))w=L(e^{z(\vare)\cdot }w)$, hence
\begin{equation*}\fl
\|L\| \ge |\vare^2z(\vare)^2-z(\vare)|\ge |\Im\,
(\vare^2z(\vare)^2-z(\vare))|=
 |\Im\, z(\vare)| |2\vare^2 \Re\, z(\vare)-1|.
\end{equation*}
Choose $\vare_0>0$ such that $ |2\vare^2 \Re\, z-1|>1/2$ for all
$z\in S,  |\vare|<\vare_0$. For $|\vare|<\vare_0$, we have
\begin{equation*}
|\Im\, z(\vare)|<2\|L\|.
\end{equation*}
Thus, for $|\vare|<\vare_0$ the solutions $z(\vare)\in S$ of $\det
\Delta_\vare (z)=0$ are necessarily inside the rectangle
$\Gamma=[\la_0-\de, \la_0 +\de_1 ]\times \big [-2\|L\|,2\|L\|\big
]$.

Now, let $F(z,\vare):=\det \Delta_\vare (z), z\in\mathbb{C},
\vare\in \mathbb{R}$. Clearly, $F(z,\vare)\to F(z,0)$ as $\vare\to
0$, for all $z\in\mathbb{C}$. Moreover, since $\Delta_\vare
(z)=\vare^2z^2I+\Delta_0(z)$, one easily deduces that the function
$F(\cdot ,\vare)$ converges uniformly to $F(\cdot , 0)$  on bounded
sets of  $\mathbb{C}$, as $\vare\to 0$.

We now apply Rouch\'e's Theorem on the boundary  $\partial\Gamma$ of
 $\Gamma$. Set $m=\min _{z\in \partial\Gamma}
|F(z,0)|>0$. For $|\vare|$ small, we have $|F(z,\vare)-F(z,0)|<m ,$
$ z\in\partial\Gamma,$ hence $F(z,\vare)$ and $F(z,0)$ have the same
number of zeros inside $\Gamma$. Thus, for $\vare>0$ sufficiently
small $\la(\vare)$ is the only solution of (\ref{e4.1}) in the strip
$S$.\ter

\smal

For (\ref{e3.2}) written as a system in $\R^{2N}$, its linearized
equation at zero is
\begin{equation}\label{e4.2}
x'(t)={\cal L}_\vare (x_t),
\end{equation}
where
\begin{equation*}
{\cal L}_\vare \pmatrix{\phi_1\cr \phi_2\cr}=\pmatrix {\phi_2(0)\cr
-{1\over {\vare^2}}L(\phi_1)+{1\over {\vare^2}} \phi_2(0)\cr},\q
\phi_1,\phi_2\in {\cal C}=C([-\tau ,0];\mathbb{R}^N).
\end{equation*}
For the linear system (\ref{e4.2}),  the characteristic equation is
given by
\begin{equation}\label{e4.3}
\det D_\vare (s)=0, \q {\rm where} \q D_\vare (s) =\pmatrix
{sI&-I\cr {1\over {\vare^2}}L(e^{s\cdot }I)& (s-{1\over
{\vare^2}})I\cr},
\end{equation}
\smal and $I$ is the $N\times N$ identity matrix. Clearly,  $\det
\Delta_\vare(s)=\vare^{2N} \det D_\vare(s)$, hence for $\vare >0$
(\ref{e4.3}) is equivalent to (\ref{e4.1}).

\begin{lem} \label{L4.2} Consider $b\in\mathbb{R}$ and $\vare_0\in (0,1)$
such that $\det \Delta_\vare(s)\ne 0$ for all $\vare \in [0,\vare_0]$ and $s$ on the vertical line
 $\Sigma =\{ s=b+iy:y\in\mathbb{R}\}$.
 Then there is $\vare_1 \in (0,\vare_0]$ such that
 $$\sup \Big \{ |s|\,  \|\Delta_\vare(s)^{-1}\|: 0\le \vare\le \vare_1, s\in \Sigma \Big \}<\infty,$$
 and
 $$\sup \Big \{ |s|\,  \|D_\vare(s)^{-1}\|: 0< \vare\le \vare_1, s\in \Sigma \Big \}<\infty.$$
\end{lem}
 {\it Proof.} For $\vare \in (0,\vare_0], s\in \Sigma$, we have
 \begin{equation}\label{e4.4}\fl
 G(\vare,s):=D_\vare(s)^{-1}=\pmatrix {\Delta_\vare(s)^{-1}&0\cr 0&\Delta_\vare(s)^{-1}\cr}
 \pmatrix {(\vare^2 s-1)I&\vare^2I\cr -L(e^{s\cdot}I)&\vare^2sI\cr}.
 \end{equation}
Clearly,
\begin{equation*}
\fl G(\vare,s)\to \pmatrix {\Delta_0(s)^{-1}&0\cr
0&\Delta_0(s)^{-1}\cr}
 \pmatrix {-I&0\cr -L(e^{s\cdot}I)&0\cr}:=G(0,s)\q {\rm as}\q \vare\to 0^+,
\end{equation*}
with $G(\vare,s)$ continuous on $[0,\vare_0]\times \Sigma.$ For
$s\in \Sigma$, we have $\|L (e^{s\cdot}I)\|\le \max (1,
e^{-b\tau})\|L\|$. It follows that
 $$
  \|G(\vare,s)\|\le  c (\vare^2|s|+1) \|\Delta_\vare(s)^{-1}\|,\q 0<\vare\le \vare_0, s\in \Sigma,$$
  for some $c>0$. Since  $\Delta_\vare(s)=(\vare^2 s^2-s)I+L(e^{s\cdot}I)$, then
  $$ \|\Delta_\vare(s)^{-1}\| \le {1\over {|s||\vare^2s-1|-\|L(e^{s\cdot}I)\|}}$$
  if $|s||\vare^2 s-1|-\|L(e^{s\cdot}I)\|>0$. Choose $\vare_1\in (0,\vare_0]$ such that $1-\vare_1^2b\ge 1/2$. Then for $s\in\Sigma$ and $0< \vare\le \vare_1$, if $|s|\ge 4\max (1,e^{-b\tau})\|L\|:=c_1$, it follows that  $|s||\vare^2 s-1|-\|L(e^{s\cdot}I)\|\ge |s|/4>0$,
thus $|s|\, \|\Delta_\vare(s)^{-1}\| \le 4$ and
$$
|s|\,\|G(\vare,s)\|\le {{c |s|(\vare^2|s|+1)}\over {|s|
|\vare^2s-1|- \|L(e^{s\cdot}I)\|}} \le {{c\vare^2|s|^2}\over {|s|
|\vare^2s-1|- \|L(e^{s\cdot}I)\|}}+4c
$$
for $\vare \in (0,\vare_1]$ and $s\in \Sigma, |s|\ge c_1.$ Now, $|s|
|\vare^2s-1|\ge |s| \sqrt {\vare^4|s|^2+1/2}\ge \vare^2 |s|^2$, and
$${{\vare^2|s|^2}\over {|s| |\vare^2s-1|- \|L(e^{s\cdot}I)\|}}\le 2$$
 if $\vare^2|s|^2\ge 2 \|L(e^{s\cdot}I)\|$; and if  $\vare^2|s|^2\le 2 \|L(e^{s\cdot}I)\|$, then
 $${{\vare^2|s|^2}\over {|s| |\vare^2s-1|- \|L(e^{s\cdot}I)\|}}\le
{{ 2 \|L(e^{s\cdot}I)\|}\over {|s|/4}}\le 2;
 $$
hence, $|s|\,\|G(\vare,s)\|\le 6c$ if $|s|\ge c_1$.

On the other hand, on  the compact set $\{(\vare, s)\in
[0,\vare_1]\times \Sigma: |s|\le c_1\}$ the continuous functions
$|s|\, \|\Delta_\vare(s)^{-1}\|$ and $|s|\, \|G(\vare,s)\|$ attain
their suprema, and the conclusion follows.\ter

We are finally in a position to prove the main result of this
section, on the existence of positive travelling wave solutions of
Eq. (\ref{e1.1}) for large wave speeds.

\begin{thm}\label{Thm4.1} Assume (H1)-(H4). Then,   there is $c^*>0$,
such that for $c>c^*$ Eq. (\ref{e1.1}) has a positive travelling
wave solution of the form $u(t,x)=\psi (ct+w\cdot x)$ for each unit
vector $w\in \mathbb{R}^p$, with $\psi(-\infty)=0,\psi(\infty)=K$.
Moreover, the components of the profile $\psi$ are increasing in the
vicinity of $-\infty$ and it satisfies $\psi(t)=O(e^{\la (\vare)t}),
\psi'(t)=O(e^{\la (\vare)t})$ at $-\infty$, where $\vare=1/c$ and
$\la (\vare)$ is the real  solution of (\ref{e4.1}) with
$\la(\vare)\to \la_0$ as $\vare\to 0^+$. \end{thm}
 {\it Proof}.
Consider Eq. (\ref{e3.2}), where $\vare=1/c$. Let $\mu \in (0,
\la_0)$ be as in the statement of Corollary \ref{Cor3.1} and satisfy
$\la_0< \mu+2\delta$ for some fixed $\delta \in (0,\mu/4)$. Suppose
also that $\la(\vare)$ is the unique solution of $\det \Delta_\vare
(z)=0$ on the strip $\mu-\delta  \le \Re\, z < 2\mu$ for all $\vare
\in [0,\vare^*]$. Fix the profiles
\begin{equation}\label{e4.5}
\psi_\vare=\psi(\vare,0), \q \psi_0(t):= u^*(t)=e^{\la_0 t}{\bf v}
+O(e^{2\mu t}), \q \vare \in [0, \vare^*],
\end{equation}
as in Corollary \ref{Cor3.3} and Theorem \ref{T2.1}. Recall that
${\bf v}$ is a positive eigenvector associated with $\la_0$. The
proof is now divided in several steps.

\smal

{\it Claim 1}.  There is $\vare_0>0$ such that  $(\psi_\vare
(t),\psi_\vare'(t))=(e^{\la (\vare)t}{\bf v}^1(\vare
),\la(\vare)e^{\la (\vare)t} {\bf v}^1(\vare ))+w_\vare (t),$ $\vare
\in [0,\vare_0)$, with continuous ${\bf v}^1(\vare ) >0, {\bf
v}^1(0) ={\bf v},$  and $w_\vare (t)=O(e^{(\la_0+\delta)t})$ at
$-\infty$. \smal

To prove the above claim, note that $x_\vare(t):=(\psi_\vare(t),
\psi_\vare '(t))$ is a solution of the system
\begin{eqnarray}\label{e4.6}
\nonumber & & x_1'(t)=x_2(t)\\
 & & \vare^2 x_2'(t)= x_2(t)-L(x_{1,t})-h_{\vare}(t),
\end{eqnarray}
where $L=Df(0)$ and
$h_{\vare}(t)=f((\psi_\vare)_t)-L((\psi_\vare)_t)$. For $\vare
>0$, equivalently we write (\ref{e4.6}) as
\begin{equation*}
x'(t)={\cal L}_\vare (x_t)-{1\over {\vare^2}}\pmatrix{0\cr
h_\vare(t)\cr},
\end{equation*}
where ${\cal L}_\vare$ is as in (\ref{e4.2}). Since $f$ is a $C^2$
function, using the Taylor formula for $f$ (cf. e.g. \cite[p.
23]{CH}), we have the estimate
$$|h_{\vare}(t)|\le \int_0^1 (1-s) \|D^2f(s (\psi_\vare)_t)\|
\|(\psi_\vare)_t\|_\infty^2\, ds, \q t\in \mathbb{R}, \vare \in
[0,\vare^*].$$ Since $\|\psi_\vare-u^*\|_\mu \to 0$ as $\vare\to
0^+$ and $u^*(t)\to 0$ as $t\to -\infty$, the continuity of $D^2f$
at 0 implies that $\| D^2f(s (\psi_\vare)_t)\|$ is uniformly bounded
on $\vare \in [0,\vare_1] \subset [0,\vare^*), \ t \leq 0, \ s \in
[0,1]$. Together with Corollary \ref{Cor3.2}, this leads to
\begin{equation}\label{e4.7}
 |x_\vare(t)|\le Ce^{\mu t},\q |h_{\vare}(t)|\le D e^{2\mu t}\q {\rm for}\q t\le 0,
\end{equation}
 for some constants $C,D$ independent of   $\vare \in [0,\vare_1]$.

We now apply Theorem \ref{TA.1}  to (\ref{e4.6}) at $-\infty$ (see
Appendix), and derive that for $\vare >0$
\begin{equation}\label{e4.8}
 x_\vare(t)=(\psi_\vare (t), \psi_\vare '(t))=z_\vare(t)+w_\vare (t)
\end{equation}
where $z_\vare(t)$ is an eigenfunction for the linear system
$x'(t)={\cal L}_\vare (x_t)$ corresponding to the set $\Lambda
_\vare=\{ z\in \mathbb{C}: \det \Delta_\vare (z)=0, \mu \le \Re\, z
<2\mu\}$ and $w_\vare(t)=O(e^{(\la_0+\delta))t})$  at $-\infty$.
From Lemma \ref{L4.1}, let $\vare$ be on some interval
$(0,\vare_0)\subset (0,\vare_1)$ such that  that $\Lambda _\vare=\{
\la(\vare)\}$. Then, $z_\vare(t)$ is an eigenfunction for
$x'(t)={\cal L}_\vare (x_t)$ associated with  the root $\la (\vare)$
of (\ref{e4.3}), hence $z_\vare(t)=e^{\la(\vare)t}{\bf v}(\vare)$
with ${\bf v}(\vare)=({\bf v}^1(\vare),{\bf v}^2(\vare))\in
\mathbb{R}^{2N}$ satisfying $D_\vare (\la(\vare)){\bf v}(\vare)=0$
for $D_\vare$ as in (\ref{e4.3}). From this we obtain ${\bf
v}^2(\vare)=\la(\vare){\bf v}^1(\vare)$ and $\Delta_\vare
(\la(\vare)){\bf v}^1(\vare)=0.$ Furthermore, from Theorem
\ref{TA.1} (with $a=\mu, b=\la_0+2\delta, \epsilon=\delta$) and
formulae (\ref{eA.5}) and (\ref{eA.7}) (adapted to the situation
$-\infty$) we get
\begin{equation*}\fl
z_\vare(t)=Res\, (e^{t\cdot}\widetilde {\bf
x}_\vare,\la(\vare))=e^{\la(\vare)t}\lim_{s\to \la(\vare)}
(s-\la(\vare))\widetilde {\bf x}_\vare(s)=e^{\la(\vare)t}({\bf
v}^1(\vare),\la(\vare){\bf v}^1(\vare))
\end{equation*}
and
\begin{equation}\label{e4.9}
w_\vare(t)={1\over {2\pi i}} \int_{b+{\delta\over
2}-i\infty}^{b+{\delta\over 2}+i\infty} e^{st}\widetilde {\bf
x}_\vare(s)\, ds,
\end{equation}
with
\begin{equation}\label{e4.10}
\widetilde {\bf x}_\vare (s):=\widetilde{ x_\vare (-\cdot
)}(-s)=G(\vare,s)\left (r_\vare(s)-\pmatrix {0\cr \widetilde {\bf
h}_\vare (s)\cr}\right ),
\end{equation}
for $G(\vare,s)=D_\vare (s)^{-1}$ and
\begin{eqnarray}\label{e4.11}
\nonumber  & & r_\vare(s)=x_\vare(0)+ {\cal L}_\vare \left (
e^{s\cdot} \int _{\cdot}^0 e^{-su} x_\vare(u)\, du \right ),\\
 & & \widetilde {\bf h}_\vare (s):=\widetilde {h_\vare (-\cdot
)}(-s)=\int_{-\infty}^0e^{-su}h_\vare(u)\, du,\q \Re\,
s<\la_0+2\delta.
\end{eqnarray}
Here,  $\widetilde{ x_\vare (-\cdot )},\widetilde{ h_\vare (-\cdot
)}$ denote  the Laplace transforms of the functions $t\mapsto
x_\vare(-t) ,t\mapsto h_\vare(-t)$, respectively. Note that
$\widetilde {\bf x}_\vare (s)$  is meromorphic for $\Re\, s<2\mu$
with a unique singularity  at $s = \la (\vare)$ which is a simple
pole  of $\Delta_\vare (s)^{-1}$  (cf. Appendix).

From (\ref{e4.8}),  we get
\begin{equation}\label{e4.12}
\psi_\vare (t)=e^{\la (\vare)t} {\bf v}^1(\vare)+w^1_\vare (t),\
\psi_\vare'(t)=\la(\vare)e^{\la (\vare)t} {\bf v}^1(\vare)+w^2_\vare
(t),
\end{equation}
with $ w^1_\vare (t)=O(e^{(\la_0+\delta) )t})$ at $-\infty$ and
$w^2_\vare (t)=(w^1_\vare (t))'$. Next, the definition of
$r_\vare(s)$ in (\ref{e4.11}) yields
$$r_\vare(s)=x_\vare(0)+\pmatrix {0\cr -{1\over\vare^2} L\Big (e^{s\cdot } \int _{\cdot}^0 e^{-st} \psi_\vare(t)\, dt\Big )\cr},$$
hence from (\ref{e4.4}) and (\ref{e4.10})   we obtain
\begin{equation}\label{e4.13}\fl
\widetilde {\bf x}_\vare (s)= G(\vare,s)x_\vare(0)- \pmatrix
{\Delta_\vare (s)^{-1} L\Big (e^{s\cdot } \int _{\cdot}^0 e^{-st}
\psi_\vare(t)\, dt\Big )\cr s\Delta_\vare (s)^{-1} L\Big (e^{s\cdot
} \int _{\cdot}^0 e^{-st} \psi_\vare(t)\, dt\Big )\cr} -\pmatrix {
\Delta_\vare (s)^{-1} \widetilde {\bf h}_\vare(s)\cr s\Delta_\vare
(s)^{-1} \widetilde {\bf h}_\vare(s)\cr}.
\end{equation}
We now extend naturally this situation for $\vare=0$.   In
(\ref{e4.5}), denote $\la(0)=\la_0, {\bf v}^1(0)={\bf v} $. Write
$\widetilde {\bf x}_\vare= (\widetilde {\bf x}_\vare^1,\widetilde
{\bf x}_\vare^2)$ for $\vare \in (0,\vare_0)$, and let  $\widetilde
{\bf x}_0^1(s)$ be defined by (\ref{e4.13})   for $\vare=0$. Note
that formula (\ref{e4.10}) can still be used to obtain $\widetilde
{\bf x}_0^1(s)$ (cf. Appendix for more details),
$$
\widetilde {\bf x}_0^1(s)=\Delta_0(s)^{-1}[r_0^1(s)+\tilde {\bf
h}_0(s)],$$ where $r_0^1(s)=u^*(0)+L\big (e^{s\cdot } \int
_{\cdot}^0 e^{-st} u^*(t)\, dt\big )=\lim_{\vare\to
0^+}r_\vare^1(s)$.

For each $\vare\in [0,\vare_0)$,  $s=\la(\vare)$ is a pole of order
one of $G(\vare,s)$, and from (\ref{e4.10}) we deduce that  for
$\vare \in [0,\vare_0)$ and  $\mu\le \Re\, s<\la_0+2\de$  the
function $A(\vare,s)$ defined by
$A(\vare,s)=(s-\la(\vare))\widetilde {\bf x}_\vare^1(s)$ for $s\ne
\la(\vare)$, $A(\vare,\la(\vare))={\bf v}^1(\vare)$ is analytic on
$s$ and  continuous on $(\vare,s)$. In particular, $\lim_{\vare\to
0^+}A(\vare,\la(\vare))= A(0,\la_0)={\bf v}>0$, hence ${\bf
v}^1(\vare)\to {\bf v}$ as $\vare\to 0^+$. Moreover,  ${\bf
v}^1(\vare)>0$ for $\vare>0$ sufficiently small. This proves Claim
1.

\med

{\it Claim 2}.  For $\vare_0^*>0$ sufficiently small, there exists a
constant $D_0>0$ such that
\begin{equation}\label{e4.14}
|w^1_\vare (t)|\le D_0e^{(\la_0+\delta)t}\q {\rm for\ all }\q t\le
0,0<\vare <\vare_0^*.
\end{equation}
\med

To prove Claim 2, once more we shall use some formulae and estimates
in the proof of Theorem \ref{TA.1} in the Appendix, changed
accordingly to account for the asymptotic behaviour at $-\infty$,
rather than $\infty$.

Define $v_\vare
(t)=(v_\vare^1(t),v_\vare^2(t))=e^{-(\la_0+\delta)t}w_\vare (t)$ and
$u_\vare
(t)=(u_\vare^1(t),u_\vare^2(t))=e^{-(\la_0+3\delta/2)t}w_\vare (t)$.
Note that $v_\vare^1(0)=w_\vare^1(0)$ an
$v_\vare^1(t)=w_\vare^1(0)-\int_t^0 (v_\vare^1)'(s)\, ds$, with
$w_\vare^1(0)=\psi_\vare(0)-{\bf v}^1(\vare)$. Then
$|w_\vare^1(0)|\le C+|{\bf v}^1(\vare)|$, where $C>0$ is as in
(\ref{e4.7}). We need to prove that $v_\vare^1(t)$ is uniformly
bounded for $t\le 0$ and $\vare>0$ small enough.
 In order to achieve this, we shall show that
there are constants $C_0,D_0>0$ and $\vare_0^*>0$, such that
\begin{equation}\label{e4.15}
\|v_\vare^1\|_{L^1(-\infty,0]}\le C_0,\q  \vare\in (0,\vare_0^*)
\end{equation}
and
\begin{equation}\label{e4.16}
\|(v_\vare^1)'\|_{L^1(-\infty,0]}\le D_0/2,\q  \vare\in
(0,\vare_0^*),
\end{equation}
so that (\ref{e4.14})  follows immediately from (\ref{e4.16})  and
$|w_\vare^1(0)|\le D_0/2$ for $\vare\in (0,\vare_0^*)$. These
uniform estimates require a careful analysis of the explicit
formulae for $w_\vare$ given in (\ref{e4.9}) and (\ref{e4.10}). We
shall prove (\ref{e4.15})  beforehand, and then use (\ref{e4.15}) to
prove (\ref{e4.16}).

First, observe that $x_\vare(t)=z_\vare(t)+w_\vare(t)$ is a solution
of  (\ref{e4.6}), with $z_\vare(t)$  being an eigenfunction for the
linear system $x'(t)={\cal L}_\vare (x_t)$, hence $w_\vare(t)$ is a
solution of system (\ref{e4.6}) as well. The definition of
$v_\vare(t)$ yields now
\begin{equation}\label{e4.17}
(v_\vare^1)'(t)=-(\la_0+\delta) v_\vare^1(t)+v_\vare^2(t)
\end{equation}
and
\begin{equation}\label{e4.18}
\vare^2(v_\vare^1)''(t)-\al \, (v_\vare^1)'(t)+P_\vare(t)=0,
\end{equation}
where $\al=1-2\vare^2(\la_0+\delta)$ and
\begin{equation}\label{e4.19}\fl
P_\vare(t)=[\vare^2(\la_0+\delta)^2-(\la_0+\delta)]
v_\vare^1(t)+L\Big (e^{(\la_0+\delta)\cdot} (v_\vare^1)_t\Big )+
e^{-(\la_0+\delta)t}h_\vare(t).
\end{equation}

Next, similarly to (\ref{eA.9}) and Remark \ref{RA1}, using
(\ref{e4.9}) and the Plancherel theorem, we obtain
$$\|v_\vare^1\|_{L^1(-\infty, 0]}=\|w^1_\vare(t)e^{-(\la_0+\delta)t}\|_{L^1(-\infty,
0]}= \| {{e^{\delta t/2}}\over {2\pi}} \int_{-\infty}^{+\infty}
e^{iut}\widetilde {\bf x}^1_\vare(\la_0 + 3\delta/2+iu)\,
du\|_{L^1(-\infty, 0]}
$$
\begin{equation}\label{e4.20}
\le C_1(\delta) \| \widetilde {\bf
x}_\vare^1(\la_0+3\delta/2-i\cdot)\|_{L^2(\R)}\q {\rm for}\q \
\vare\in (0,\vare_0).
\end{equation}
From (\ref{e4.7}), $ |x_\vare(0)|$  is uniformly bounded for $\vare
\in (0,\vare_0)$; and since $2\mu>\la_0+2\delta$, (\ref{e4.7}) also
implies that $\| L\Big (e^{s\cdot } \int _{\cdot}^0 e^{-st}
\psi_\vare(t)\, dt\Big )\|$ and $|\widetilde {\bf h}_\vare(s)|$ are
uniformly bounded for $\vare \in (0,\vare_0)$ and $\Re\,
s=\la_0+3\delta/2$. Using now Lemma \ref{L4.2}, from (\ref{e4.13})
we conclude that there is $K_1>0$ such that
\begin{equation*}
|\widetilde {\bf x}_\vare^1(s)|\le K_1/|s|,\q {\rm for }\
s=b+{{3\delta}\over 2}+iy,\ y\in\mathbb{R},\ {\rm and}\
0<\vare<\vare_0,
\end{equation*}
from which we derive
$$\| \widetilde {\bf x}_\vare^1(\la_0+3\delta/2-i\cdot)\|_{L^2(\R)}^2\le
\int_{\R} {K_1\over { (b+{{3\delta}\over 2})^2+y^2}}\, dy<\infty.$$
Together with (\ref{e4.20}), this leads to the estimate
(\ref{e4.15}) with $\vare_0^*= \vare_0$.

In order to prove (\ref{e4.16}), we now use (\ref{e4.18})  and  some
ideas from Aguerrea et al. \cite[Lemma 4.1]{ATV}. After partial
integration of Eq. (\ref{e4.18}), we obtain
\begin{equation}\label{e4.21}
(v_\vare^1)'(t)=e^{(\al/\vare^2)t} (v_\vare^1)'(0)+
{{e^{(\al/\vare^2)t}}\over \vare^2} \int_t^0 e^{-(\al/\vare^2)s}
P_\vare (s)\, ds.
\end{equation}
Let $\vare>0$ be small, so that $\al>0$. From (\ref{e4.17})  and
Claim 1,
$$ (v_\vare^1)'(0)=-(\la_0+\de) w_\vare^1(0)+w_\vare^2(0)=
-(\la_0+\de)\psi_\vare(0)+(\la_0+\de -\la(\vare)){\bf
v}^1(\vare)+\psi_\vare'(0),$$ with ${\bf v}^1(\vare)\to {\bf v}$ as
$\vare\to 0^+$.   Corollary \ref{Cor3.2} implies that
$|(v_\vare^1)'(0)|$ is uniformly bounded on $\vare >0$ sufficiently
small, hence
\begin{equation}\label{e4.22}
|(v_\vare^1)'(0)|\int_{-\infty}^0 e^{(\al/\vare^2)t}\, dt\le K_2,
\end{equation}
for some $K_2>0$ and $\vare >0$ sufficiently small. On the other
hand, interchanging the order of  integration leads to
$$
{1\over {\vare^2}}\int_\sigma^0 e^{(\al/\vare^2)t} \left ( \int_t^0
e^{-(\al/\vare^2)s}P_\vare (s)\, ds\right ) dt={1\over {\vare^2}}
\int_\sigma^0 e^{-(\al/\vare^2)s}P_\vare (s)\left ( \int_\sigma^s
e^{(\al/\vare^2)t}\, dt\right ) ds, $$ hence
$$ {1\over {\vare^2}}\int_{-\infty}^0 e^{(\al/\vare^2)t} \left ( \int_t^0 e^{-(\al/\vare^2)s}|P_\vare (s)|\, ds\right ) dt\le {1\over \al} \int _{-\infty}^0 |P_\vare (s)|\, ds.$$
From  (\ref{e4.7}), (\ref{e4.15})  and the definition of $P_\vare
(s)$ in (\ref{e4.19}) , we easily see that
$$\|P_\vare\|_{L^1(-\infty,0]}\le K_3,\q 0<\vare <\vare_0^*$$
for some $\vare_0^*>0$. Together with (\ref{e4.21}) and
(\ref{e4.22}), this yields the estimate (\ref{e4.16}), and therefore
Claim 2 is proven.

\med

We finally prove:

\med

{\it Claim 3}. There is $\vare_1^*>0$ such that $\psi_\vare (t)>0$
for $t\in\mathbb{R}$ and $\vare \in (0,\vare_1^*)$.

\med

From Claims 1 and 2, for $\vare >0$ small enough
$$\psi_\vare (t)\ge e^{\la(\vare)t}{\bf v}^1(\vare)-D_0e^{(\la_0+\delta)t} \vec{\bf 1}
\ge e^{\la(\vare)t}[{\bf v}^1(\vare)-D_0 e^{(\delta /2)t}\vec {\bf
1}],\q t\le 0,$$ where $\vec {\bf 1}=(1,\dots, 1)$ and
$|\la(\vare)-\la_0|<\delta/2$. Choose $T^*\le 0$ and $\vare_1^*>0$
such that
$${\bf v}^1(\vare)> {\bf v}/2 \q {\rm if}\q 0<\vare<\vare_1^*, \q {\rm and}\q D_0 e^{(\delta /2)T^*}\vec {\bf 1}\le {\bf v}/4.$$
Then
$$\psi_\vare (t)\ge e^{\la(\vare)t} {\bf v}/4>0, \q 0<\vare<\vare_1^*,\ t\le T^*.$$
On the other hand, since $\| \psi_\vare -u^*\|_\infty\to 0$ as
$\vare\to 0^+$, we define $\eta:=\inf_{t\ge T^*} u^*(t)>0$, and
suppose that $\vare_1^*$ was chosen so that
 $\| \psi_\vare -u^*\|_\infty<\eta$ for $0<\vare<\vare_1^*$.
 It follows that $\psi_\vare (t)>0$ for all $t\in\mathbb{R}$ and $\vare\in (0,\vare_1^*)$. The proof of
  the theorem is complete.\ter

   \med

The above proof shows that the requirement that there is a positive
eigenvector ${\bf v}\in \mathbb{R}^N$ for the dominant
characteristic value $\la_0$ of (\ref{e2.1}) is crucial to deduce
the {\it positiveness} of the travelling wave fronts, for large wave
speeds. Nevertheless, the {\it existence} of such waves and their
asymptotic behaviour at $-\infty$ can be deduced from our Theorem
\ref{Thm3.2}, as well as the auxiliary results in Section 3,  and
the proof of the above Claim 1. We summarize these remarks in the
following theorem:
\begin{thm} \label{Thm4.2} Assume (H1), (H2)(i),  and\vskip 0cm
  (i) for Eq. (\ref{e1.2}) the equilibrium $u=K$ is localy asymptotically stable;\vskip 0cm
  (ii) for Eq. (\ref{e1.2}), the linearized equation about 0 has a real characteristic root $\la_0>0$, which is simple and dominant;\vskip 0cm
  (iii) Eq. (\ref{e1.2}) has a heteroclinic solution $u^*(t), t\in\mathbb{R}$,
with $u^*(-\infty)=0, u(\infty)=K$ and $u^*(t)= O(e^{\la_0 t})$ at
$-\infty$.\vskip 0cm Then,  there is $c^*>0$, such that for $c>c^*$,
Eq. (\ref{e1.1}) has a travelling wave solution of the form
$u(t,x)=\psi (ct+w\cdot x)$ for each unit vector $w\in
\mathbb{R}^p$, with
 $\psi(t)=O(e^{\la (\vare)t}), \psi'(t)=O(e^{\la (\vare)t})$ at $-\infty$, where $\vare=1/c$
 and $\la (\vare)$ is the real  solution of (\ref{e4.1}) with  $\la(\vare)\to \la_0$ as $\vare\to 0^+$.
\end{thm}

\section{Applications}
\subsection{A diffusive generalized logistic equation with
distributed delay}
 As a first  application, we consider a scalar reaction-diffusion equation with distributed delays in the reaction-terms, which includes
 the Fisher-KPP equation with delay as a particular case.

 Let ${\cal C}=C([-\tau,0];\mathbb{R}), \tau >0,$ and consider
\begin{equation}\label{e5.1}
 {{\p u}\over {\p t}}(t,x)={{\p u^2}\over {\p x^2}} (t,x)+bu(t,x)[1-L(u_t(\cdot ,x))],\q t\in \mathbb{R},\ x\in \mathbb{R},
\end{equation}
where $b>0$ and $L:{\cal C}\to \mathbb{R}$ is a nonzero positive
linear operator, i.e., $L\ne 0$ is linear and $L(\var)\ge 0$
whenever $\var \ge 0$.  In particular, $L$ is  bounded and
$\|L\|=L(1)>0$.

The corresponding delayed ODE model,
\begin{equation}\label{e5.2}
u'(t)=bu(t)[1-Lu_t],\q t\in\mathbb{R},
\end{equation}
has two equilibria, $u=0$ and $u=L(1)^{-1}:=K$. Here $f$ in (1.2)
reads as $f(\var)=b\var (0)[1-L(\var)], \var\in {\cal C}$, for which
it is easy to verify that conditions (H1) and (H2)  are satisfied.
The linearized equation about zero is the ODE $u'(t)=bu(t)$, with
characteristic value $b>0$.

Now, $u'(t)=-bKLu_t$ is the linearized equation about $K$, with
characteristic equation $P(\la ):=\la +bK L(e^{\la \cdot})=0$. If
$\la=i\omega, \omega>0,$ is a solution of $P(\la )=0$, then $L(\cos
(\omega \cdot))=0$ and $0= \omega +bK L(\sin (\omega \cdot))\ge
\omega-b$. If $\omega\tau <\pi /2$, we deduce that there is $\de >0$
such that  $\cos (\omega \th)\ge \de$ for $\th \in [-\tau,0]$, hence
$L(\cos (\omega \cdot))\ge \de L(1)>0$, which is a contradiction. On
the other hand,  if $b\tau\le 3/2$, from  \cite{F} we conclude that
the positive  equilibrium $K$ is globally attractive in the set of
all positive solutions of (\ref{e5.2}) with initial conditions
$\var\in{\cal C}_+,\var (0)>0$. We thus conclude that (H3) holds if
$b\tau\le 3/2$. Therefore, the following result is  an immediate
consequence of Theorem \ref{Thm4.1}
\begin{thm}\label{T5.1} If $b\tau\le 3/2$, there exists $c^*>0$ such that
for $c>c^*$  equation (\ref{e5.1}) has a positive travelling wave
solution $u(t,x)=\psi (x+ct)$ with $\psi (-\infty )=0,\psi
(\infty)=K$. Moreover, $\psi$ is increasing in a vicinity of
$-\infty$ and it has the asymptotic decay $\psi (t)=O(e^{b(c)t}),
\psi '(t)=O(e^{b(c)t})$ at $-\infty$, where
$b(c)=2bc/(c+\sqrt{c^2-4b})$. \end{thm}

We note that (\ref{e5.1}) includes as a particular case the
Fisher-KPP equation with a single  delay,
\begin{equation}\label{e5.3}
{{\p u}\over {\p t}}(t,x)={{\p u^2}\over {\p x^2}}
(t,x)+bu(t,x)[1-u(t-\tau ,x)/K],\q t\in \mathbb{R},\ x\in
\mathbb{R}.
\end{equation}
By using a pair of upper-lower solutions and a monotone iterative
method, Wu and Zou \cite{WZ} proved that if $c>2\sqrt b$, then there
exists $\tau^*(c)>0$ such that, for a delay $\tau\le\tau ^*(c)$,
(\ref{e5.3}) has a non-decreasing travelling wave front connecting 0
to $K$ with wave speed $c$. Our approach does not allow us to
determine the minimal wave speed $c^*$, but, on the contrary, we
explicitly exhibit the maximal delay $\tau^*=3/(2b)$, under which we
can assure the existence of such positive (but not necessarily
monotone) travelling solutions.

Travelling waves for the Fisher-KPP equation (\ref{e5.3}) with $b=1$
were also considered in \cite[Corollary 6.6]{FHW}, where it was
shown that, for $\tau \le e^{-1}$,  there exists  $c^*>0$ such that
there is a travelling wave front with wave speed $c>c^*$. We remark
that Theorem \ref{T5.1} applied to (\ref{e5.3}) with $b=1$ clearly
improves this result, since it guarantees the existence of
travelling wave solutions for $\tau\le 3/2$, and, most relevant in
biological terms, it does assert that such travelling waves are {\it
positive}.

For some other recent results and references on the Fisher-KPP
equation, see \cite{BNPR, GT}.

\subsection{A chemostat model with delayed  growth response }

Consider the following model for the growth of bacteria in a
well-stirred chemostat supplied by a single essential nutrient (cf.
Ellermeyer \cite{E1} and Ellermeyer et al. \cite{E2}):
\begin{eqnarray}\label{e5.4}
\nonumber  & & S'(t)=D(S^0-S(t))-f(S(t))u(t),\\
 & & u'(t)=e^{-D\tau} f(S(t-\tau))u(t-\tau)-Du(t).
\end{eqnarray}
Here $S(t)$ and $u(t)$ are the concentration of nutrient in the
growth vessel and the biomass concentration of bacteria at time $t$,
respectively, $D>0$ is the dilution rate of the chemostat, $S^0>0$
is the input concentration of nutrient, and $\tau\ge 0$ is the delay
in the growth response, to account for the lag in the nutrient
conversion into biomass due to  cellular absorption; $f$ is the
specific functional response for the bacteria, and typically the
Michaelis-Menten response is chosen, $ f(s)=ms/(a+s),\ s\ge 0, $
with $m,a>0$. More generally, one can consider a continuously
differentiable and bounded function $f:[0,\infty)\to [0,\infty)$
with
\begin{equation}\label{e5.5}
f(0)=0,\q f'(s)>0\ {\rm for}\ s\in\mathbb{R}.
\end{equation}
For an unstirred chemostat, the nutrient is added to the vessel but
not mixed, so one has to introduce diffusion terms. The diffusion
rates  $d_1,d_2>0$ for the nutrient and the organisms may be
different, and model (\ref{e5.4}) becomes
\begin{eqnarray}\label{e5.6}
\nonumber  & & {{\p S}\over {\p t}}(t,x)=d_1\Delta
S(t,x)+D(S^0-S(t,x))-f(S(t,x))u(t,x)\\
 & & {{\p u}\over {\p t}}(t,x)=d_2\Delta u(t,x)+e^{-D\tau}
f(S(t-\tau,x))u(t-\tau,x)-Du(t,x),
\end{eqnarray}
for $t\in \mathbb{R}$ and $x\in (0,L)\ (L>0)$ (or more generally
$x\in\Omega$, where $\Omega \subset \mathbb{R}^3$ is an open
domain). There is an extensive literature on  chemostat models with
``delayed growth response",  with and without diffusion. We refer to
\cite{BR,E1,E2,GK,SW,WXR},  for results, other related chemostat
models, biological explanations, and further references.

For both (\ref{e5.4}) and (\ref{e5.6}), there is always  the
equilibrium $(S^0,0)$, corresponding to the ``washout state". If
\begin{equation}\label{e5.7}
f(S^0)>De^{D\tau},
\end{equation}
there is another nonnegative equilibrium $(\bar S,\bar u)$, called
the ``survival state", given by $(\bar S,\bar
u)=(f^{-1}(De^{D\tau}), e^{-D\tau}(S^0-\bar S))$. Condition
(\ref{e5.7}) imposes a restriction on the size of the time-delay
$\tau$, which should satisfy $\tau<D^{-1}\log(f(S^0)D^{-1})$  for
$f$ as in (\ref{e5.5}). Moreover,  (\ref{e5.7}) implies that  the
equilibrium $(S^0,0)$ of (\ref{e5.4}) is  unstable and $(\bar S,\bar
u)$ is asymptotically stable and a global attractor of all solutions
with initial conditions $(S_0,u_0)=(\phi_1,\phi_2)\in {\cal C}_+$,
$\phi_2(0)>0$ \cite{E1,E2}.

In order to apply our results, we first observe that both the
positive cone ${\cal C}_+$ and the set $\{ (\phi_1,\phi_2)\in {\cal
C}_+: \phi_2(0)> 0, \phi_1(0)< S^0\}$ are positively invariant for
(\ref{e5.4}). Translating the washout state to the origin, by
setting $s(t)=S^0-S(t)$, we rewrite (\ref{e5.4}) as
\begin{eqnarray}\label{e5.8}
\nonumber  & & s'(t)=-Ds(t)+f(S^0-s(t))u(t)\\
 & & u'(t)=e^{-D\tau} f(S^0-s(t))u(t-\tau)-Du(t).
\end{eqnarray}
If $f$ is $C^2$-smooth and (\ref{e5.7}) holds, then Eq. (\ref{e5.8})
satisfies (H1) and (H2). The set $A=\{ (\phi_1,\phi_2)\in {\cal
C}_+: \phi_2(0)> 0, 0<\phi_1(0)\le S^0\}$ is positively invariant,
and (\ref{e5.8}) has equilibria $E_0=(0,0)$ and $K:= (\bar s,\bar
u)>0$, with $\bar s=S^0-\bar S, \bar u=e^{-D\tau} \bar s$, the first
one being unstable and the second one being locally stable and a
global attractor of all solutions  with initial conditions in  $A$
(cf. \cite{E1}). It remains to verify that (H4) holds.

The characteristic equation for the linearization of (\ref{e5.8}) at
$(0,0)$ is
\begin{equation}\label{e5.9}
(\la+D)(\la+D-e^{-D\tau} f(S^0)e^{-\la \tau})=0.
\end{equation}
Define $h(x)=(x+D)e^{(x+D)\tau}, x\in \mathbb{R}$. Under
(\ref{e5.7}), there are two real roots of (\ref{e5.9}), $-D$ and
$\la_0$, where $\la_0>0$ is the unique solution of $h(x)=f(S^0)$.
For $\la\notin \mathbb{R}$ a solution of (\ref{e5.9}), we have
$h(\Re\, \la)<f(S^0)$. Since $h'(x)>0$ for $x>0$, if follows that
$\Re\, \la <\la_0$, and we conclude that $\la_0$ is a dominant
eigenvalue for the linearization of (\ref{e5.8}) about $(0,0)$, with
${\bf v}=(f(S^0)(\la_0+D)^{-1}, 1)>0$ as an associated eigenvector.
From Theorems \ref{T2.1} and \ref{Thm4.1}, we therefore obtain the
following result:

\begin{thm}\label{Thm5.2} Consider Eq. (\ref{e5.4}), where $S^0,D,\tau >0$,
the function  $f:[0,\infty)\to [0,\infty)$ is bounded,  $C^2$-smooth
and satisfies (\ref{e5.5}) and (\ref{e5.7}). Then, there exists a
heteroclinic solution $(S^*(t), u^*(t))$ connecting the washout
state $(S^0,0)$ to the survival state $(\bar S,\bar u)$ of
(\ref{e5.4}), with $0<S^*(t)<S^0, u^*(t)>0$ for $t\in\mathbb{R}$ and
$(S^*(t), u^*(t))=(S^0,0)+O(e^{\la_0 t})$ at $-\infty$, where
$\la_0>0$ satisfies $\la_0+D=f(S^0)e^{-(\la_0+D)\tau}$. For the
diffusion model (\ref{e5.6}) with $d_1,d_2>0$, there is $c^*>0$ such
that for $c>c^*$ (\ref{e5.8})   has a positive travelling wave
solution of the form $(S(t,x),u(t,x))=(\psi_1(ct+x),\psi_2(ct+x))$,
with $\psi_1(-\infty)=S^0,\psi_2(-\infty)=0$ and
$\psi_1(\infty)=\bar S,\psi_2(\infty)=\bar u$; moreover,
$\psi_1'(t)<0,\psi_2'(t)>0$ in the vicinity of $-\infty$, and
$(\psi_1(t),\psi_2(t))=(S^0,0)+O(e^{z (c) t})$ at $-\infty$, where
$z(c)$ is the real solution of
\begin{equation}\label{e5.10}\fl
z^2-cz -c(D-e^{-(D+\la)\tau}f(S^0))=0.
\end{equation}
\end{thm}
{\it Proof.} Let $\vare=1/c$, for $c>0$ large. With the  notation in
(\ref{e4.3}),  we have that
\begin{equation}\label{e5.11}\fl
\det D_\vare (\la)={1\over {\vare^2 d_1d_2}}( \vare d_1\la^2-\la -D)
\Big (\vare \la^2-\la -D+ e^{-(D+\la)\tau}f(S^0)\Big ).
\end{equation}
Define $\la(\vare)$ as the real solution of (\ref{e5.11}) such that
$\la(\vare)\to \la_0$ as $\vare\to 0^+$. Then $z(c)=\la(\vare)$
satisfies (\ref{e5.10}).\ter \smal

Theorem \ref{Thm5.2}  asserts  the existence of  {\it positive}
travelling wavefronts for (\ref{e5.6}) with $S(t,x)<S^0$ for all
$t\in\R,x\in\R$. Here, due to the change of variables $s=S-s^0$, the
positivity of the component $s(t)$ (or $s(t,x)$) translates  as the
nutrient concentration being smaller than $S^0$. We emphasize that
biologically significant solutions of (\ref{e5.4}) and (\ref{e5.6})
must be positive and have a nutrient concentration $S$ not larger
than the input concentration $S^0$.

\begin{rem}\label{R5.1} The existence of a positive eigenvector ${\bf v}$
associated with the dominant eigenvalue $\la_0$, as prescribed in
(H4), may seem a quite restrictive requirement, since it is not
satisfied by many populations dynamics systems, namely Kolmogorov
type models with $N>1$. We however observe that   if the
characteristic matrix  for (\ref{e2.1}) at $\la_0$,
$\Delta_0(\la_0)$, is an irreducible matrix with non-negative
off-diagonal entries, then there is a positive eigenvector for
$\Delta_0(\la_0)$ associated with $\la_0$ (see e.g. \cite[p.
258]{SW}). This property will be exploited in a forthcoming paper,
where Theorems \ref{Thm3.2} and \ref{Thm4.1} will be applied to
several population models.
\end{rem}

\section{Appendix}

In this appendix, we extend Proposition 7.1 of Mallet-Paret
\cite{MP} to systems with distributed delays.

Consider the  FDE
\begin{equation}\label{eA.1}
x'(t)=L_0x_t+h(t),\q t\in\mathbb{R}
\end{equation}
and the homogeneous system
\begin{equation}\label{eA.2}
x'(t)=L_0x_t,
\end{equation}
where $L_0:C([-\tau, 0];\mathbb{R}^N)\to \mathbb{R}^N$ is a bounded
linear operator and $h:\mathbb{R}\to\mathbb{R}^N$ is continuous. For
(\ref{eA.2}), write the characteristic equation
\begin{equation*}
\det \Delta_0(s)=0,\q {\rm where}\q \Delta_0(s)=sI-L_0(e^{s\cdot
}I).
\end{equation*}
It is well known that the solutions of the characteristic equation
are exactly the eigenvalues for the homogeneous system (\ref{eA.2}),
i.e., the eigenvalues for the infinitesimal generator $A$ associated
with the semiflow of (\ref{eA.2}). Furthermore, the spectrum $\sigma
(A)$ of $A$ is only composed of the point spectrum.

\begin{lem}\label{LA.1} If $f$ is a holomorphic function on a disc $\{
s: |s-\la|<\epsilon\} $, where $\la$ is an eigenvalue of
(\ref{eA.2}) and $\epsilon>0$ is small, then
\begin{equation*}
x(t)=Res\,  (e^{t\cdot} \Delta_0^{-1}f,\la)={1\over {2\pi
i}}\int_{|s-\la|=\epsilon} e^{st} \Delta_0(s)^{-1}f(s)\, ds
\end{equation*}
is an eigenfunction of (\ref{eA.2}) corresponding to $\la$.
\end{lem}
{\it Proof}. The proof follows the arguments  of  Mallet-Paret
\cite[Section 7]{MP}, so we omit it.\ter

\begin{thm}\label{TA.1} Let $x(t)$ be a solution of (\ref{eA.1}) on
$[T,\infty)$ for some $T\in\mathbb{R}$. Assume there are
$a,b\in\mathbb{R}, a<b$, such that
\begin{equation*}
x(t)=O(e^{-at}),\q h(t)=O(e^{-bt})\q {\rm as}\q t\to\infty.
\end{equation*}
Then, for every $\epsilon >0$, we have
\begin{equation*}x(t)=z(t)+O(e^{-(b-\epsilon)t})\q {\rm as}\q t\to\infty,
\end{equation*}
where $z(t)$ is an eigenfunction of (\ref{eA.2}) associated with the
set of eigenvalues $\Lambda=\{ \la \in \sigma (A):-b<\Re\, \la \le
-a\}.$ Analogously, if $x(t)$ is a solution of (\ref{eA.1}) on
$(-\infty, T]$ for some $T\in\mathbb{R}$ and
\begin{equation*}
x(t)=O(e^{at}),\q h(t)=O(e^{bt})\q {\rm as}\q t\to -\infty,
\end{equation*}
with $a<b$, then for every $\epsilon >0$ we have
\begin{equation*}
x(t)=z(t)+O(e^{(b-\epsilon)t})\q {\rm as}\q t\to -\infty,
\end{equation*}
where $z(t)$ is an eigenfunction of (\ref{eA.2}) associated with the
set of eigenvalues $\Lambda=\{ \la \in \sigma (A):a\le \Re\, \la
<b\}.$
\end{thm}
{\it Proof}. We only prove the result for $+\infty$; for $-\infty$
it is analogous. Without loss of generality, take $T=0$. In what
follows, for $f:[0,\infty)\to \mathbb{C}, f(t)=O(e^{-at})$ at
$+\infty$, we denote the Laplace transform of $f$ by
$$({\cal L}f)(s)=\tilde f(s)=\int_0^\infty e^{-st} f(t)\, dt,\q {\rm for}\q \Re\, s>-a.$$

Write $L_0(\var)=\int_{-\tau}^0 d\eta (\th) \var (\th)$, where $\eta
(\th)$ is an $N\times N$ matrix-valued  function of bounded
variation.   Applying the Laplace transform to (\ref{eA.1}), we
obtain
\begin{equation}\label{eA.3}
-x(0)+s\tilde x(s)={\cal L} (L_0 x_t)(s)+\tilde h(s), \q \Re\, s>-a,
\end{equation}
where \begin{eqnarray}\label{eA.4}
\nonumber\fl  & & {\cal L} (L_0 x_t)(s)=\int _0^\infty e^{-st} \left
(\int _{-\tau}^0 d\eta (\th) x(t+\th )\right ) dt=\int _{-\tau}^0
d\eta (\th)
\left ( \int _0^\infty e^{-s(t+\th)} x(t+\th)\, dt\right )\\
 \fl & & = \int _{-\tau}^0 d\eta (\th)  \left ( \int_\th ^0 e^{-su}
x(u)\, du +\tilde x(s)\right )=L_0\left ( e^{s\cdot} \int _{\cdot}^0
e^{-su} x(u)\, du \right )+L_0(e^{s\cdot} I)\tilde x(s).
\end{eqnarray}
From $(\ref{eA.3}),(\ref{eA.4})$, we obtain
$$
 \Delta_0(s)\tilde x(s)=r(s)+\tilde h(s),\q {\rm where}\q
r(s)=x(0)+L_0\left ( e^{s\cdot} \int _{\cdot}^0 e^{-su} x(u)\, du
\right).
$$
We observe that $r(s)$ is an entire function, $\tilde h(s)$ is
defined and holomorphic for  $\Re\, s>-b$, and
\begin{equation}\label{eA.5}
\tilde x(s)=\Delta_0(s)^{-1} [r(s)+\tilde h(s)]
\end{equation}
is holomorphic for $\Re\, s>-b$ with the exception of finitely many
poles.

Take $\epsilon >0$ and $k>-a$.  On any strip of the form
$-b+\epsilon\le \Re\, \la \le k$, the functions $r(s)$ and $\tilde
h(s)$  are uniformly bounded. Since $k$ is greater than the real
part of all singularities of $\tilde x(s)$, we can use the inverse
formula for the Laplace transform,
\begin{equation}\label{eA.6}
x(t)={1\over {2\pi i}} \int_{k-i\infty}^{k+i\infty} e^{st} \tilde
x(s)\, ds.
\end{equation}
Choose $\epsilon>0$ small such that $-b+\epsilon/2<k$ and $\sigma
(A)\cap \{ s:-b<\Re\, s\le -b+\epsilon/2\} =\emptyset.$ Note that in
the strip $-b\le \Re\, s\le -a$, the only possible poles of $\tilde
x$ lie on $-b+\epsilon/2<\Re\, s\le -a$.

In the strip $-b+\epsilon/2\le \Re\, s\le k$, the functions
$|r(s)|,| \tilde h(s)|$ are bounded, and $\|L_0(e^{s\cdot }I)\|\le
\max (1, e^{-k\tau})\|L_0\|$, hence the  operator norm for the
inverse of $\Delta_0(s)=sI-L_0(e^{s\cdot }I)$ satisfies
\begin{equation*}
\| \Delta_0(s)^{-1}\| \le {1\over {|s|-\max (1,
e^{-k\tau})\|L_0\|}}\q {\rm for}\q |s|>\max (1, e^{-k\tau})\|L_0\|.
\end{equation*}
From this estimate and (\ref{eA.5}), we conclude that $| e^{st} \tilde x(s)| \to 0$ as $|\Im\, s|\to \infty$, uniformly
  in the strip $-b+\epsilon/2\le \Re\, s\le k$, and that $\tilde x(s)$ is $L^2$-integrable in any straight line
$s=x_0+iy,y\in\mathbb{R}$, for any fixed $x_0\in [ -b+\epsilon/2,
k]$.
 We may shift the path of integration in (\ref{eA.6}) to the left, and obtain
 \begin{equation*}
x(t)=z(t)+w(t),\q {\rm where}
 \end{equation*}
 \begin{equation}\label{eA.7}
z(t)=\sum_{\la \in\Lambda} Res\, (e^{t\cdot} \tilde x,\la),\q
 w(t)={1\over {2\pi i}} \int_{-b+{\epsilon\over 2}-i\infty}^{-b+{\epsilon\over 2}+i\infty} e^{st} \tilde x(s)\, ds.
 \end{equation}
From the previous lemma, $z(t)$ is an eigenfunction of (\ref{eA.2})
associated with the set of eigenvalues in $\Lambda$. It remains to
prove that
\begin{equation}\label{eA.8}
w(t)=O(e^{-(b-\epsilon)t})\q {\rm at}\q +\infty.
\end{equation}

 Define $u(t)=e^{(b-\epsilon/2)t}w(t), v(t)=e^{(b-\epsilon)t}w(t)$.
 We first prove that $v\in L^2[0,\infty)$. Here,  $L^p[0,\infty)=L^p([0,\infty);\mathbb{C}^N),p=1,2$.
 We have
\begin{equation*}
\fl u(t)={1\over {2\pi i}} \int_{-b+{\epsilon\over
2}-i\infty}^{-b+{\epsilon\over 2}
 +i\infty} e^{(s+b-\epsilon/2)t} \tilde x(s)\, ds={1\over {2\pi}}\int _{\R} e^{-ist} \tilde x(-b+\epsilon/2-is)\, ds.
\end{equation*}
By Plancherel Theorem, $u(t)\in L^2[0,\infty)$, so that
$v(t)=e^{-\epsilon t/2} u(t)$ implies
 \begin{equation}\label{eA.9}\fl
 \|v\|_{L^1[0,\infty)}
 \le \|u\|_{L^2[0,\infty)}
 \|e^{-\epsilon t/2}\|_{L^2[0,\infty)}
 \le C\| \tilde x(-b+\epsilon/2+i\cdot )\|_{L^2(\mathbb{R})} <\infty,
 \end{equation}
for some $C>0$. Hence, $v\in L^1[0,\infty).$ Define now
\begin{equation*}
V(t)=L_0(e^{-(b-\epsilon)\cdot }v_t)=\int_{-\tau}^0d\eta(\th)
e^{-(b-\epsilon)\th}v(t+\th).
\end{equation*}
Then,
\begin{equation}\label{eA.10}
 \int_0^{+\infty} |V(t)|\, dt\le \max (1, e^{(b-\epsilon)\tau})\|L_0\|\|v\|_{L^1[-\tau,\infty)}.
\end{equation}
 In particular, $V\in L^1[0,\infty)$.
We now observe that $x(t)=z(t)+w(t)$, with $x(t)$ a solution of
(\ref{eA.1}) and $z(t)$ a solution of (\ref{eA.2}). Hence $w(t)$
satisfies (\ref{eA.1}), $w'(t)=L_0(w_t)+h(t)$, and we obtain
\begin{equation*}
v'(t)=(b-\epsilon)v(t)+L_0(e^{-(b-\epsilon)\cdot}v_t)+e^{(b-\epsilon)t}h(t),
\end{equation*}
with $e^{(b-\epsilon)t}h(t)=O(e^{-\epsilon t})$. We conclude that
$v'\in L^1[0,\infty)$. Since $|v(t)|\le |v(0)|+\int_0^t |v'(s)|\,
ds$ for $ t\ge 0$, then $v$ is bounded on $[0,\infty)$, and
(\ref{eA.8}) holds.\ter

\begin{rem}\label{RA1} For the situation $x(t)=O(e^{-at}), h(t)=O(e^{-bt})\ (a<b)$ as
$t\to\infty$, denote $v(t)=(v_1(t),\dots,v_N(t))$  as in the above
proof. Clearly,  one can obtain componentwise estimates similar to
(\ref{eA.9}) or (\ref{eA.10}). In fact, one  concludes that for
$\epsilon>0$ small such that $\sigma (A)\cap \{ s:-b<\Re\, s\le
-b+\epsilon/2\} =\emptyset$ and $t\ge 0, j=1,\dots,N,$
\begin{equation*}
\|v_j\|_{L^1[0,\infty)}\le  (2\pi \sqrt{\epsilon})^{-1}\| \tilde
x_j(-b+\epsilon/2+i\cdot )\|_{L^2(\mathbb{R})}
\end{equation*}
and $|v_j(t)|\le |v_j(0)|+\| v_j'\|_{L^1[0,\infty)}$ with
\begin{equation*}
\| v_j'\|_{L^1[0,\infty)}\le C  \| \tilde
x_j(-b+\epsilon/2+i\cdot)\|_{L^2(\mathbb{R})}+\|e^{(b-\epsilon)\cdot}h_j\|_{L^1[0,\infty)},
\end{equation*}
where $C={{(b-\epsilon)+e^{(b-\epsilon/) \tau}\|L_0\|}\over {2\pi
\sqrt {\epsilon}}} $. Similar estimates hold for the case
$x(t)=O(e^{at}), h(t)=O(e^{bt})$ at $ -\infty.$ \end{rem}

\ack This research was supported by  FCT (Portugal), Financiamento
Base 2008-ISFL-1-209 (Teresa Faria) and by FONDECYT (Chile),
projects 7080045 (Teresa Faria) and 1071053 (Sergei Trofimchuk). S.
Trofimchuk was also partially supported by CONICYT (Chile) through
PBCT program ACT-56 and by the University of Talca, program
``Reticulados y Ecuaciones".

\end{document}